\documentclass[aop,citesort,MSNbibl,dvips]{arximspdf}
\usepackage{graphics}
\doi{10.1214/09-AOP504}
\volume{38}
\issue{3}
\pubyear{2010}
\firstpage{1019}
\lastpage{1061}

\makeatletter

\newtheorem{theo}{Theorem}[section]
\newtheorem{lmm}{Lemma}[section]
\newproclaim{df}{Definition}[section]
\newtheorem{prop}{Proposition}[section]
\newproclaim{rem}{Remark}[section]
\newproclaim{Condition}{Condition}

\makeatother

\begin{document}
\begin{frontmatter}

\title{Quenched invariance principle for the Knudsen stochastic
billiard in a random tube}
\runtitle{Stochastic billiards in a random tube}

\begin{aug}
\author[A]{\fnms{Francis} \snm{Comets}\thanksref{t1}\corref{}\ead[label=e1]{comets@math.jussieu.fr}},
\author[B]{\fnms{Serguei} \snm{Popov}\thanksref{t2,t5}\ead[label=e2]{popov@ime.unicamp.br}},
\author[C]{\fnms{Gunter M.} \snm{Sch\"{u}tz}\thanksref{t3}\ead[label=e3]{g.schuetz@fz-juelich.de}}\\ and
\author[B]{\fnms{Marina} \snm{Vachkovskaia}\thanksref{t4,t5}\ead[label=e4]{marinav@ime.unicamp.br}}
\runauthor{Comets, Popov, Sch\"{u}tz and Vachkovskaia}
\affiliation{Universit{\'{e}} Paris Diderot (Paris 7),
University of Campinas--UNICAMP,
Forschungszentrum J\"{u}lich GmbH,
Institut f\"{u}r Festk\"{o}rperforschung
and University of Campinas--UNICAMP}
\address[A]{F. Comets\\
Universit{\'{e}} Paris Diderot (Paris 7)\\
UFR de Math\'{e}matiques\\
case 7012 (Site Chevaleret)\\
F-75205
%2 place Jussieu\\
Paris Cedex 13\\
France\\
\printead{e1}} %adresu isvedimo komanda gale!
\address[B]{S. Popov\\
M. Vachkovskaia\\
University of Campinas--UNICAMP\\
\printead{e2}\\
\phantom{E-mail: }\printead*{e4}}
\address[C]{G. M. Sch\"{u}tz\\
Forschungszentrum J\"{u}lich GmbH\\
Institut f\"{u}r Festk\"{o}rperforschung\\
\printead{e3}}
\end{aug}

\thankstext{t1}{Supported in part by CNRS (UMR 7599
``Probabilit{\'{e}}s et Mod{\`{e}}les Al{\'{e}}atoires'')
and ANR (grant POLINTBIO).}

\thankstext{t2}{Supported in part by CNPq (300328/2005--2).}

\thankstext{t3}{Supported by DFG (Schu 827/5--2, Priority programme SPP
1155).}

\thankstext{t4}{Supported in part by CNPq (304561/2006--1).}

\thankstext{t5}{Supported by CNPq (471925/2006--3), FAPESP
(04/07276--2) and CAPES/DAAD
(Probral).}

\pdfauthor{Francis Comets, Serguei Popov, Gunter M. Schutz, Marina Vachkovskaia}

% HISTORY:
\received{\smonth{11} \syear{2008}}
\revised{\smonth{10} \syear{2009}}

% ABSTRACT
%
\begin{abstract}
We consider a stochastic billiard in a random tube which stretches to
infinity in the direction of the first coordinate. This random tube is
stationary and ergodic,
and also it is supposed to be in some sense
well behaved.
The stochastic billiard can be described as follows: when strictly
inside the tube, the particle moves
straight with constant speed. Upon hitting the boundary,
it is reflected randomly, according to the
cosine law: the density of the outgoing direction is proportional
to the cosine of the angle between
this direction and the normal vector. We also consider the
discrete-time random walk formed
by the particle's positions at the moments of hitting the boundary.
Under the condition of existence of the second moment
of the projected jump length with respect to the stationary
measure for the environment seen from the particle,
we prove the quenched invariance principles for the projected
trajectories of the random walk and the stochastic billiard.
\end{abstract}

% KEYWORDS
%

\end{frontmatter}

%s1 ###
\section{Introduction}
\label{s_intro}
The so-called Knudsen regime in gas dynamics describes a very
dilute gas confined between solid walls. The gas is dilute in the
sense that the mean free path of gas molecules, that is, the typical
distance travelled between collisions of gas molecules, is much
larger than the typical distance between consecutive collisions of
the gas molecules with the walls. Hence molecules interact
predominantly with the walls, and the interaction among themselves
can be neglected.
A typical setting where this is relevant is in the motion of absorbed
guest molecules in the pores of a microporous solid
where both the pore
diameter and the typical number of guest molecules
inside the pores are small.

On molecular scale the wall-molecule interaction is usually
rather complicated and very difficult to handle explicitly.
Hence one
resorts to a stochastic description in which gas molecules, from now
on referred to as particles, move ballistically between collisions
with the walls where they interact in a random fashion. In the
Knudsen model one assumes that particles are pointlike and that the
kinetic energy of a particle is conserved in a collision with the
wall, but its direction of motion changes randomly.
The law of this random reflection is taken to be the cosine law
where the outgoing direction is cosine-distributed relative to the
surface normal at the point where the particle hits the wall.
For a motivation of this choice, sometimes also called Lambert
reflection law, see \cite{FY}. Notice that this dynamic implies
that the incoming direction is not relevant
and is ``forgotten'' once a collision has happened.
Thus this process defines a Markov chain which we call
``Knudsen stochastic billiard'' (KSB).
The random sequence of hitting points is referred to as
``Knudsen random walk'' (KRW) \cite{CPSV1}.

Pores in microporous solids may have a very complicated surface.
Among the many possibilities, an elongated tube-shaped pore surface
has recently
attracted very considerable attention \cite{ZRBCK}.
A three-dimensional connected
network of ``tubes'' may be regarded as constituting the entire
(nonsimply connected) interior empty space of microporous grain
in which particles can move. In this setting parts of the surface of
the individual tubes are open and connect to a neighboring pore so
that particles can move from one to another pore. It
is of great interest to study the large scale motion of a molecule
along the direction in which the tube
has its greatest elongation. We
think of the direction of longest elongation of a single tube as the
first coordinate in $d$-dimensional Euclidean space. Together with
the locally, usually very complicated surface of pores, this
leads us to introduce the notion of a random tube with a random
surface to be defined precisely below.

Knudsen motion in a tube has been studied heuristically for simple
regular geometries such as a circular pipe and also numerically for
self-similar random surfaces (see, e.g., \cite{CD,CM,RZBK}). For the
straight infinite pipe, it is not difficult to see that
in dimensions larger than two, the mean-square displacement grows
asymptotically linearly in time, that is, diffusively, while in two
dimensions the motion is superdiffusive due to sufficiently large
probability for very long flights between collisions.
Interestingly though, rigorous work on this conceptually simple
problem is rare. In fact, it is not even
established under which conditions on the pore surface the motion
of a Knudsen particle has diffusive mean square displacement and
converges to Brownian motion. Indeed, it is probable that
one may construct counterexamples to Brownian motion in three or
more dimensions without having to
invent physically pathological pore
surfaces, but considering a nonstationary (expanding or shrinking)
tube instead. We refer here to the work \cite{MVW} (there, only the
two-dimensional case is considered, but it seems reasonable that one
may expect similar phenomena in the general case as well). Even for
the stationary tube, it seems reasonable that the presence of
bottlenecks with random (not bounded away from $0$) width may cause
the process to be subdiffusive.

In this work we shall define a class of single
infinite random tubes
for which we prove convergence of KSB and KRW to Brownian motion.
Together with related earlier and on-going work, this lays the
foundation for addressing more subtle issues that arise in the study
of motion in open domains which are finite and which allow for
injection and extraction
of particles.
This latter problem is of
great importance for studying the relation between the
so-called self-diffusion coefficient, given by the asymptotic
mean square displacement of a particle in an infinite tube
under equilibrium conditions, on the one hand, and the
transport diffusion coefficient on the other hand. The
transport diffusion coefficient is given by the
nonequilibrium flux in a finite open tube with a density
gradient between two open ends of the tube. Often only
one of these quantities can be measured experimentally,
but both may be required for the correct interpretation of
other experimental data. Hence one would like to
investigate whether both diffusion coefficients are equal
and specify conditions under which this the case. This
is the subject of a forthcoming paper \cite{CPSV3}. Here we focus
on the question of diffusion in the infinite tube in
the absence of fluxes.

For a description of our strategy we first come to the modeling of
the infinite tube. The tube stretches to infinity
in the first direction. The collection of its sections in the
transverse direction can be thought as a ``well-behaved'' function
$\omega\dvtx \alpha\mapsto\omega_\alpha$ of the first
coordinate $\alpha$ which values are subsets of ${\mathbb R}^{d-1}$.
We assume
that the boundary of the tube is Lipschitz and almost everywhere
differentiable, in order to define the process.
We also assume the transverse sections are bounded and that there
exist
no long dead ends or too thin bottlenecks.
For our long-time asymptotics of the walk, we need some large-scale
homogeneity assumptions on the tube; it is quite natural to
assume that the process $\omega=( \omega_\alpha;
\alpha\in{\mathbb R})$ is random, stationary and ergodic.
Now, the process essentially looks like a one-dimensional
random walk in a random environment, defined by
random conductances since Knudsen random walk is reversible.
The tube serves as a random environment for our walk. The tube, as
seen from the walker, is a Markov process which has a (reversible)
invariant law.
From this picture, we understand that the random medium is
homogenized on large scales,
and that, for almost every environment, the walk is asymptotically
Gaussian in the longitudinal direction. More precisely, we will
prove that after the usual rescaling the trajectory of the KRW
converges weakly to Brownian motion with some diffusion constant
(and from this we deduce also the analogous result for the KSB).
This will be done by showing ergodicity for the environment
seen from the particle and using the classical ``corrector
approach'' adapted from \cite{Koz}.

The point of view of the particle has become useful
\cite{KV,DFGW,BS} in the study of reversible random walks
in a random environment and in obtaining the central limit theorem for
the annealed law (i.e., in the mean with respect to the
environment). The authors from \cite{FSS} obtain an (annealed)
invariance principle for a random walk on a random point process in
the Euclidean space, yielding an upper bound on the effective
diffusivity which agrees with the predictions of Mott's law.
The corrector approach, by correcting the walk into a martingale in
a fixed environment, has been widely used to obtain the quenched
invariance principle for the simple random walk on the infinite
cluster in a supercritical percolation \cite{SS} in dimension
$d \geq4$, \cite{BB} for $d \geq2$.
These last two references apply to walks defined by random
conductances which are
bounded and bounded away from 0. The authors in \cite{BP}
and \cite{M} gave a proof under the only
condition of bounded conductances.

In this paper we will leave untouched the questions
of deriving heat kernel estimates or spectral estimates (see \cite{Bar}
and \cite{FM}, respectively) for the corresponding results
for the simple random walk on the infinite cluster.
We will not need such estimates to show that the corrector
can be neglected in the limit. Instead we will benefit from the
essentially one-dimensional structure
of our problem and use the ergodic theorem;
this last ingredient will require
introducing reference points in Section \ref{s_refpoints}
below to recover stationarity.

The paper is organized as follows. In Section \ref{s_notations}
we formally define the random tube and construct the stochastic
billiard and the random walk and then formulate our results.
In Section \ref{s_environment} the process of the environment seen
from the particle is defined and its properties are discussed.
Namely, in Section \ref{s_env_discrete} we define
the process of the environment seen from the discrete-time random
walk and then prove that this process is reversible with an
explicit reversible measure. For the sake of completeness,
we do the same for the continuous-time stochastic billiard in
Section~\ref{s_env_continuous}, even though the results of this
section are not needed for the subsequent discussion.
In Section \ref{s_constr_corrector} we construct the corrector
function, and in Section \ref{s_refpoints} we show that the corrector
behaves sublinearly along a certain stationary sequence of reference
points and that one may control the fluctuations of the corrector
outside this sequence.
Based on the machinery developed in Section \ref{s_environment},
we give the proofs of our results in Section \ref{s_proofs}. In
Section \ref{s_proof_inv_pr} we prove the quenched
invariance principle for the discrete-time random walk, and in
Section \ref{s_b_finite} we discuss the question of
finiteness of the averaged second moment of the projected jump
length. In Section~\ref{s_proof_inv_pr_cont} we prove the quenched
invariance principle for the continuous-time stochastic billiard,
also obtaining an explicit relation between the corresponding
diffusion constants.
Finally, in the \hyperref[app]{Appendix} we discuss the general case of random tubes
where vertical walls are allowed.

%s2 ###
\section{Notation and results}
\label{s_notations}
Let us formally define the random tube in ${\mathbb R}^d$, \mbox{$d\geq2$}.
In this paper, ${\mathbb R}^{d-1}$ will always stand for the linear
subspace of ${\mathbb R}^d$ which is perpendicular to the first coordinate
vector $e$; we use the notation \mbox{$\|\cdot\|$} for the Euclidean norm
in ${\mathbb R}^d$ or ${\mathbb R}^{d-1}$. For $k\in\{d-1,d\}$
let ${\mathcal B}(x,\varepsilon)=\{y\in{\mathbb R}^k\dvtx\|x-y\|
<\varepsilon\}$ be the open
$\varepsilon$-neighborhood of $x\in{\mathbb R}^k$. Define
${\mathbb S}^{d-1}=\{y\in{\mathbb R}^d\dvtx\|y\|=1\}$ to be the unit
sphere in ${\mathbb R}^d$,
and let ${\mathbb S}^{d-2}={\mathbb S}^{d-1}\cap{\mathbb R}^{d-1}$ be
the unit sphere
in ${\mathbb R}^{d-1}$. We write $|A|$ for the $k$-dimensional Lebesgue
measure in case $A\subset{\mathbb R}^k$ and $k$-dimensional Hausdorff
measure in case $A\subset{\mathbb S}^k$. Let
\[
{\mathbb S}_h = \{w\in{\mathbb S}^{d-1}\dvtx h\cdot w > 0\}
\]
be the half-sphere looking in the direction $h$.
For $x\in{\mathbb R}^d$, sometimes it will be convenient to write
$x=(\alpha,u)$; $\alpha$ being the first coordinate of $x$ and
$u\in{\mathbb R}^{d-1}$; then $\alpha=x\cdot e$,
and we write $u={\mathcal U}x$; ${\mathcal U}$ being the projector on
${\mathbb R}^{d-1}$.
Fix some positive constant ${\widehat M}$, and define
%e1 ###
%
\begin{equation}
\label{def_Lambda}
\Lambda= \{u\in{\mathbb R}^{d-1} \dvtx \|u\| \leq{\widehat M}\}.
\end{equation}

Let $A$ be an open
%FC: connected
domain in ${\mathbb R}^{d-1}$ or ${\mathbb R}^d$.
We denote by $\partial A$ the boundary of $A$ and
by $\bar A = A\cup\partial A$ the closure of $A$.
\begin{df}
\label{def_Lipschitz}
Let $k\in\{d-1,d\}$, and suppose that $A$ is an open
%FC: connected
domain in ${\mathbb R}^k$. We say that $\partial A$ is
$({\hat\varepsilon},{\hat L})$-Lipschitz, if for any $x\in\partial A$
there exist an affine isometry ${\mathfrak I}_x \dvtx {\mathbb R}^k\to
{\mathbb R}^k$ and a
function $f_x\dvtx{\mathbb R}^{k-1}\to{\mathbb R}$ such that:
\begin{itemize}
\item$f_x$ satisfies Lipschitz condition with constant ${\hat L}$,
i.e., $|f_x(z)-f_x(z')| < {\hat L}\|z-z'\|$ for all $z,z'$;
\item${\mathfrak I}_x x = 0$, $f_x(0)=0$ and
\[
{\mathfrak I}_x\bigl(A\cap{\mathcal B}(x,{\hat\varepsilon})\bigr) = \bigl\{z\in
{\mathcal B}(0,{\hat\varepsilon})\dvtx
z^{(k)} > f_x\bigl(z^{(1)},\ldots,z^{(k-1)}\bigr)\bigr\}.
\]
\end{itemize}
In the degenerate case $k=1$ we adopt the convention
that $\partial A$ is $({\hat\varepsilon},{\hat L})$-Lipschitz
%FC: for any positive ${\hat\eps},{\hat L}$.
for any positive ${\hat L}$ if points in $\partial A$ have a mutual
distance larger than ${\hat\varepsilon}$.
\end{df}

%FC:
Fix ${\widehat M}>0$, and define ${\mathfrak E}_n$ for $n \geq1$
to be the set of all open domains $A$ such that $A\subset\Lambda$
and $\partial A$ is $(1/n,n)$-Lipschitz.
We turn ${\mathfrak E}= \bigcup_{n \geq1}{\mathfrak E}_n$ into a
metric space by defining
the distance between $A$ and $B$ to be equal to
$|(A\setminus B)\cup(B\setminus A)|$, making ${\mathfrak E}$ a Polish space.
Let $\Omega={\mathcal C}_{\mathfrak E}({\mathbb R})$ be the space of
all continuous
functions ${\mathbb R}\mapsto{\mathfrak E}$; let ${\mathcal A}$ be the
sigma-algebra
generated by the cylinder sets with respect to the Borel
sigma-algebra on ${\mathfrak E}$, and let ${\mathbb P}$ be a
probability measure on
$(\Omega,{\mathcal A})$. This defines a ${\mathfrak E}$-valued process
$\omega=(\omega_\alpha, \alpha\in{\mathbb R})$ with continuous
trajectories.
Write $\theta_\alpha$ for the spatial shift: $\theta_\alpha
\omega_{\cdot} = \omega_{\cdot+\alpha}$. We suppose that the
process $\omega$ is stationary and ergodic with respect to the
family of shifts $(\theta_\alpha,\alpha\in{\mathbb R})$.
With a slight abuse of notation, we denote also by
\[
\omega= \{(\alpha,u)\in{\mathbb R}^d \dvtx u\in\omega_\alpha\}
\]
the random domain (``tube'') where the billiard lives.
Intuitively, $\omega_\alpha$ is the ``slice'' obtained
by crossing $\omega$ with the hyperplane
$\{\alpha\}\times{\mathbb R}^{d-1}$. One can check that, under
Condition \ref{ConditionL}
below, the domain $\omega$ is an open subset of
${\mathbb R}^d$, and we will assume that it is connected.

We assume the following:
\renewcommand{\theCondition}{L}
\begin{Condition}\label{ConditionL}
There exist ${\tilde\varepsilon},{\tilde L}$ such that $\partial
\omega$
is $({\tilde\varepsilon},{\tilde L})$-Lipschitz (in the sense of
Definition \ref{def_Lipschitz}) ${\mathbb P}$-a.s.
\end{Condition}

Let $\mu^{\omega}_\alpha$ be the $(d-2)$-dimensional Hausdorff
measure on
$\partial\omega_\alpha$
(in the case $d=2$, $\mu^{\omega}_\alpha$ is simply the counting measure),
and denote $\mu^{\omega}_{\alpha,\beta} = \mu^{\omega}_\alpha
\otimes\mu^{\omega}_\beta$.
Since it always holds that $\partial\omega_\alpha\subset\Lambda$, we
can regard $\mu^{\omega}_\alpha$ as a measure on $\Lambda$ (with
$\operatorname{supp}
\mu^{\omega}_\alpha= \partial\omega_\alpha$), and $\mu^{\omega
}_{\alpha,\beta}$
as a measure on $\Lambda^2$ (with $\operatorname{supp}\mu^{\omega
}_{\alpha,\beta} =
\partial\omega_\alpha\times\partial\omega_\beta$). Also, we denote
by $\nu^{\omega}$ the $(d-1)$-dimensional Hausdorff measure
on $\partial\omega$; from Condition \ref{ConditionL} one obtains that $\nu^{\omega}$
is locally finite.

We keep the usual notation $dx, dv, dh, \ldots$ for the
$(d-1)$-dimensional Lebesgue measure
on $\Lambda$ (usually restricted to $\omega_\alpha$ for
some $\alpha$) or the
surface measure on ${\mathbb S}^{d-1}$.

%f1 ###
%
\begin{figure}

\includegraphics{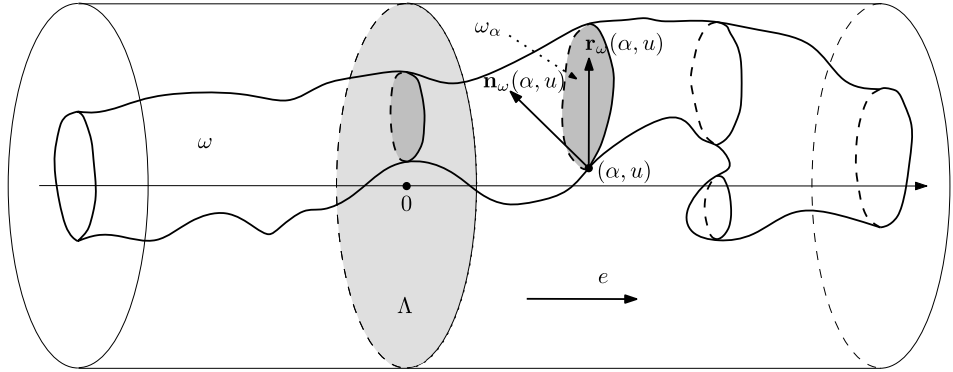}

\caption{Notation for the random tube.
Note that the sections may be disconnected.}
\label{f_tube}
\end{figure}

\begin{figure}

\includegraphics{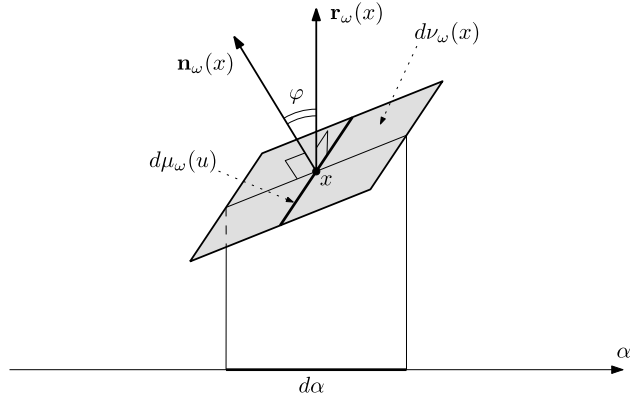}

\caption{On formula (\protect\ref{differentials});
note that $\kappa_{\alpha,u}=\cos\phi$.}
\label{f_differentials}
\end{figure}

For all $x=(\alpha,u)\in\partial\omega$ where they exist,
define the normal
vector ${\mathbf n}_{\omega}(x)={\mathbf n}_{\omega}(\alpha,u)\in
{\mathbb S}^{d-1}$ pointing inside the
domain $\omega$
and the vector ${\mathbf r}_{\omega}(x)={\mathbf r}_{\omega}(\alpha
,u)\in{\mathbb S}^{d-2}$ which is the
normal vector at $u\in\partial\omega_\alpha$ pointing inside the
domain $\omega_\alpha$
(in fact, ${\mathbf r}_{\omega}$ is the normalized projection of
${\mathbf n}_{\omega}$
onto ${\mathbb R}^{d-1}$) (see Figure \ref{f_tube}).
Denote also
\[
\kappa_x=\kappa_{\alpha,u}={\mathbf n}_{\omega}(x)\cdot{\mathbf
r}_{\omega}(x).
\]
Observe that $\kappa$ also depends on $\omega$, but we prefer not
to write it explicitly in order not to overload the notation.
Define the set of regular points
\[
{\mathcal R}_{\omega}= \{x\in\partial\omega\dvtx \partial\omega\mbox
{ is
continuously differentiable in }x, |{\mathbf n}_{\omega}(x)\cdot
e|\neq1\}.
\]
Note that $\kappa_x\in(0,1]$ for all $x\in{\mathcal R}_{\omega}$.
Clearly, it holds that
%e2 ###
%
\begin{equation}
\label{differentials}
d\nu^{\omega}(\alpha,u) = \kappa_{\alpha,u}^{-1} \,d\mu^{\omega
}_\alpha(u)\,
d\alpha
\end{equation}
(see Figure \ref{f_differentials}).

%f2 ###
%

We suppose that the following condition holds:
\renewcommand{\theCondition}{R}
\begin{Condition}\label{ConditionR}
We have $\nu^{\omega}(\partial\omega\setminus{\mathcal R}_{\omega
})=0$, ${\mathbb P}$-a.s.
\end{Condition}

We say that $y\in\bar\omega$ is \textit{seen from} $x\in\bar\omega$
if there exists $h\in{\mathbb S}^{d-1}$ and $t_0>0$ such that
$x+th\in\omega$ for all $t\in(0,t_0)$ and $x+t_0 h = y$.
Clearly, if $y$ is seen from $x$, then $x$ is seen from $y$,
and we write ``$x \stackrel{\omega}{\leftrightarrow}y$'' when this occurs.

One of the main objects of study in this paper is the
Knudsen random walk (KRW) $\xi=(\xi_n)_{n \in{\mathbb N}}$ which is a
discrete time Markov process on $\partial\omega$
(cf. \cite{CPSV1}). It is defined through its transition density $K$:
for $x,y\in\partial\omega$,
%e3 ###
%
\begin{equation}
\label{def_K}
K(x,y) = \gamma_d \frac{ ((y-x)\cdot{\mathbf n}_{\omega}(x) )
((x-y)\cdot{\mathbf n}_{\omega}(y) )}{\|x-y\|^{d+1}}
{\mathbb I}{\{x,y\in{\mathcal R}_{\omega}, x \stackrel{\omega
}{\leftrightarrow}y\}},
\end{equation}
where $\gamma_d = (\int_{{\mathbb S}_{e}} h\cdot e \,dh )^{-1}$
is the normalizing constant.
This means that, being ${\mathtt P}_{\omega},{\mathtt E}_{\omega}$
the quenched (i.e.,
with fixed $\omega$) probability and expectation, for any $x\in
{\mathcal R}_{\omega}$
and any measurable $B\subset\partial\omega$
we have
\[
{\mathtt P}_{\omega}[\xi_{n+1}\in B \mid\xi_n=x] = \int_B K(x,y)\,
d\nu^{\omega}(y).
\]
Following \cite{CPSV1}, we shortly explain why this Markov chain
is of natural interest. From $\xi_{n}=x$, the next step
$\xi_{n+1}=y$ is performed by picking randomly the direction
$h=(y-x)/\|y-x\|$ of the step
according to Knudsen's cosine density
$\gamma_d h \cdot{\mathbf n}_{\omega}(x) \,dh$ on the half
unit-sphere looking toward the interior of the domain. By
elementary\vspace*{1pt} geometric considerations, one can check that
$d\nu^{\omega}(y)= (h \cdot{\mathbf n}_{\omega}(y))^{-1} \|y-x\|
^{d-1}\,dh$ and recover the
previous formulas.
%FC: We also refer to the Knudsen random walk as the
%random walk with cosine
%reflection law, since it can be easily seen
%from \eqref{def_K} that the
%density of the outgoing direction is proportional
%to the cosine of the angle
%between this
%direction and the normal vector (e.g., formula (4)
%in \cite{CPSV1}).

Let us define also
%e4 ###
%
\begin{equation}
\label{K->Phi}
\Phi(\alpha,u,\beta,v)= (\kappa_{\alpha,u}\kappa_{\beta,v})^{-1}
K ((\alpha,u),(\beta,v) ).
\end{equation}

%FC: As observed in \cite{CPSV1},
From (\ref{def_K}) we see that
$K(\cdot,\cdot)$ is symmetric, that is, $K(x,y)=K(y,x)$ for all
$x,y\in{\mathcal R}_{\omega}$;
consequently, $\Phi$ has this property as well:
%e5 ###
%
\begin{equation}
\label{symmetry_Phi}
\Phi(\alpha,u,\beta,v) = \Phi(\beta,v,\alpha,u)
\qquad\mbox{for all }\alpha,\beta\in{\mathbb R},
u\in\partial\omega_\alpha, v\in\partial\omega_\beta.
\end{equation}
Clearly, both $K$ and $\Phi$ depend on $\omega$ as well,
but we usually do not
indicate this in the notation in order to keep them simple.
When we have to do it,
we write $K^\omega, \Phi^\omega$ instead of $K,\Phi$.
For any $\gamma$ we have
%e6 ###
%
\begin{equation}
\label{shift_Phi}
\Phi^{\theta_\gamma\omega}(\alpha,u,\beta,v) =
\Phi^\omega(\alpha+\gamma,u, \beta+\gamma,v).
\end{equation}
Moreover, the symmetry implies that
%e8 ###
%e7 ###
%
\begin{eqnarray}
\label{symm_K}
K^\omega((0,u),y ) &=& K^{\theta_{y\cdot e}\omega}
\bigl((0,{\mathcal U}y),(-y\cdot e,u) \bigr),
\\
\label{symm_Phi}
\Phi^\omega(0,u,\alpha,v) &=& \Phi^{\theta_\alpha\omega}
(0,v,-\alpha,u).
\end{eqnarray}

We need also to assume the following technical condition:
\renewcommand{\theCondition}{P}
\begin{Condition}\label{ConditionP}
There exist constants $N,\varepsilon,\delta$ such that for ${\mathbb
P}$-almost
every $\omega$, for any $x,y\in{\mathcal R}_{\omega}$ with
$|(x-y)\cdot e|\leq2$
there exist $B_1,\ldots,B_n\subset\partial\omega$,
$n\leq N-1$ with $\nu^{\omega}(B_i)\geq\delta$ for all $i=1,\ldots
,n$ and
such that:
\begin{itemize}
\item$K(x,z)\geq\varepsilon$ for all $z\in B_1$,
\item$K(y,z)\geq\varepsilon$ for all $z\in B_n$,
\item$K(z,z')\geq\varepsilon$ for all $z\in B_i$, $z'\in B_{i+1}$,
$i=1,\ldots,n-1$
\end{itemize}
[if $N=1$ we only require that $K(x,y)\geq\varepsilon$].
In other words, there exists a ``thick'' path of length at
most $N$ joining $x$ and $y$.
\end{Condition}

Now, following \cite{CPSV1}, we define also the
Knudsen stochastic billiard (KSB)
$(X,V)$. First, we do that for the process starting on the
boundary $\partial\omega$
%FC: from the point $x_0\in\partial\omega$.
from the point $x_0\in{\mathcal R}_{\omega}\subset\partial\omega$.
Let $x_0=\xi_0,\xi_1,\xi_2,\xi_3,\ldots$
be the trajectory of the random walk, and define
\[
\tau_n = \sum_{k=1}^n \|\xi_k-\xi_{k-1}\|.
\]
Then, for $t\in[\tau_n,\tau_{n+1})$, define
\[
X_t=\xi_n+(\xi_{n+1}-\xi_n)\frac{t-\tau_n}{\|\xi_{n+1}-\xi_n\|}.
\]
The quantity $X_t$ stands for the position of the particle
at time $t$.
Since $(X_t)_{t \geq0}$ is not a Markov process by itself,
we define also the c\`{a}dl\`{a}g version
of the motion direction at time $t$,
\[
V_t = \lim_{\varepsilon\downarrow0}\frac{X_{t+\varepsilon
}-X_t}{\varepsilon}.
\]
Then, $V_t\in{\mathbb S}^{d-1}$ and the couple $(X_t,V_t)_{t \geq0}$
is a Markov process. Of course, we can define also the stochastic
billiard starting from the interior of $\omega$ by specifying its
initial position $X_0$ and initial direction $V_0$.

Define
\[
{\mathfrak S}= \{(\omega,u) \dvtx \omega\in\Omega, u\in\partial\omega
_0\}.
\]
One of the most important objects in this paper is the
probability measure ${\mathbb Q}$ on ${\mathfrak S}$ defined by
%e9 ###
%
\begin{equation}
\label{def_Q}
d{\mathbb Q}(\omega,u) = \frac{1}{{\mathcal Z}} \kappa_{0,u}^{-1}
\,d\mu^{\omega}_{0}(u)
\,d{\mathbb P}(\omega),
\end{equation}
where ${\mathcal Z}=\int_\Omega d{\mathbb P}
\int_\Lambda\kappa_{0,u}^{-1}\,d\mu^{\omega}_0(u)$
is the normalizing constant. (We will show that ${\mathbb Q}$ is the
invariant law of the environment seen from the walker.)
To see that ${\mathcal Z}$ is finite, note that
${\mathcal Z}= \int_\Omega d{\mathbb P}
\int_{0}^1 d\alpha\int_\Lambda
\kappa_{\alpha,u}^{-1}\,d\mu^{\omega}_\alpha(u)$
by translation invariance, that is,
the expected surface area of the tube
restricted to the slab $[0,1] \times{\mathbb R}^{d-1}$
which is finite by Condition \ref{ConditionL}.
Let $L^2({\mathfrak S})$ be the space of ${\mathbb Q}$-square
integrable functions
$f\dvtx {\mathfrak S}\mapsto{\mathbb R}$.
We use the notation $\langle f \rangle_{{\mathbb Q}}$ for the
${\mathbb Q}
$-expectation of $f$:
\[
%%\label{def_expect}
\langle f \rangle_{{\mathbb Q}} = \frac{1}{{\mathcal Z}}\int_{\Omega
}d{\mathbb P}\int
_{\Lambda}
d\mu^{\omega}_0(u)\, \kappa_{0,u}^{-1}
f(\omega,u)
\]
and we
define the scalar product $\langle\cdot, \cdot\rangle_{{\mathbb
Q}}$ in
$L^2({\mathfrak S})$ by
%e10 ###
%
\begin{equation}
\label{def_scalar}
\langle f, g \rangle_{{\mathbb Q}} = \frac{1}{{\mathcal Z}}\int
_{\Omega}d{\mathbb P}
\int_{\Lambda}
d\mu^{\omega}_0(u)\, \kappa_{0,u}^{-1}
f(\omega,u) g(\omega,u).
\end{equation}
Note that $\langle f \rangle_{{\mathbb Q}} = \langle{\mathbf1}, f
\rangle
_{{\mathbb Q}}$ where ${\mathbf
1}(\omega,u)=1$ for all $\omega,u$.

Now, for $(\beta,u)\in{\mathcal R}_{\omega}$ we define the local
drift and the second
moment of the jump projected on the horizontal direction:
%e11 ###
%
\begin{eqnarray}\label{def_drift}
\Delta_\beta(\omega,u) &=& {\mathtt E}_{\omega}\bigl((\xi_1-\xi_0)\cdot
e \mid
\xi_0=(\beta,u) \bigr)\nonumber\\
&=& \int_{\partial\omega}(x\cdot e - \beta) K((\beta,u), x)
\,d\nu^{\omega}(x)
%% \label{def_drift_K}
\\
&=& \int_{-\infty}^{+\infty} (\alpha-\beta)\, d\alpha
\int_\Lambda d\mu^{\omega}_\alpha(v)\,
\kappa_{\beta,u} \Phi(\beta,u,\alpha,v),
\nonumber\\
b_\beta(\omega,u) &=& {\mathtt E}_{\omega}\bigl(\bigl((\xi_1-\xi_0)\cdot e\bigr)^2
\mid
\xi_0=(\beta,u) \bigr)\nonumber\\
&=& \int_{\partial\omega}(x\cdot e - \beta)^2 K((\beta,u), x)
\,d\nu^{\omega}(x)
\nonumber%% \label{def_second_moment_K}
\\
&=& \int_{-\infty}^{+\infty} (\alpha-\beta)^2 \,d\alpha
\int_\Lambda d\mu^{\omega}_\alpha(v)\,
\kappa_{\beta,u} \Phi(\beta,u,\alpha,v).
\nonumber%% \label{def_second_moment}
\end{eqnarray}
When $\beta=0$, we write simply $\Delta(\omega,u)$ and
$b(\omega,u)$
instead of $\Delta_0(\omega,u)$ and $b_0(\omega,u)$.
In Section \ref{s_environment} we show that $\langle\Delta\rangle
_{{\mathbb Q}}=0$
[see (\ref{g_Delta})].

Let $Z^{(m)}_\cdot$ be the polygonal interpolation of
$n/m \mapsto m^{-1/2}\xi_n\cdot e$.
Our main result is the quenched invariance principle
for the horizontal projection of the random walk.
\begin{theo}
\label{t_q_invar_princ}
Assume Conditions \ref{ConditionL}, \ref{ConditionP}, \ref{ConditionR}, and suppose that
%e12 ###
%
\begin{equation}
\label{finite_2nd_moment}
\langle b \rangle_{{\mathbb Q}} < \infty.
\end{equation}
Then, there exists a constant $\sigma>0$ such that for ${\mathbb P}$-almost
all $\omega$, for any starting point from ${\mathcal R}_{\omega}$,
$\sigma^{-1}Z^{(m)}_\cdot$
converges in law, under ${\mathtt P}_{\omega}$,
to Brownian motion as \mbox{$m\to\infty$}.
\end{theo}

The constant $\sigma$ is defined by (\ref{def:sigma}) below.
Next, we obtain the corresponding result for the continuous time
Knudsen stochastic billiard.
Define ${\hat Z}_t^{(s)} = s^{-1/2}X_{st}\cdot e$.
Recall also a notation from \cite{CPSV1}:
for $x\in\omega$, $v\in{\mathbb S}^{d-1}$,
define (with the convention $\inf\varnothing=\infty$)
\[
{\mathsf h}_x(v) = x+v\inf\{t>0 \dvtx x+tv \in\partial\omega\}
\in\{\partial\omega,\infty\}
\]
so that ${\mathsf h}_x(v)$ is the next point where the particle
hits the boundary when starting at the location $x$
with the direction $v$.
\begin{theo}
\label{t_q_invar_princ_cont}
Assume Conditions \ref{ConditionL}, \ref{ConditionP}, \ref{ConditionR}, and suppose
that (\ref{finite_2nd_moment}) holds.
Denote
\[
{\hat\sigma} = \frac{\sigma\Gamma({d}/{2}+1){\mathcal Z}}
{\pi^{1/2}\Gamma(({d+1})/{2})d}
\biggl(\int_\Omega|\omega_0| \,d{\mathbb P}\biggr)^{-1},
\]
where $\sigma$ is from Theorem \ref{t_q_invar_princ}.
Then, for ${\mathbb P}$-almost all $\omega$,
for any initial conditions $(x_0,v_0)$ such that
${\mathsf h}_{x_0}(v_0)\in{\mathcal R}_{\omega}$,
${\hat\sigma}^{-1}{\hat Z}^{(s)}_\cdot$
converges in law to Brownian motion as $s\to\infty$.
\end{theo}

Next, we discuss the question of validity
of (\ref{finite_2nd_moment}).
\begin{prop}
\label{pr_suff_2nd_moment}
If $d\geq3$ then (\ref{finite_2nd_moment}) holds.
\end{prop}

If $d=2$, then one cannot expect (\ref{finite_2nd_moment})
to be valid in the situation when $\omega$ contains an infinite
straight cylinder. Indeed, we have the following:
\begin{prop}
\label{pr_infinite_d=2}
In the two-dimensional case, suppose that there exists an interval
$S\subset\Lambda$ such that ${\mathbb R}\times S \subset\omega$ for
${\mathbb P}$-a.a. $\omega$. Then $\langle b \rangle_{{\mathbb
Q}}=\infty$.
\end{prop}

% In fact,
% in this situation it was analytically argued
% that the random walk is
% superdiffusive.
%FC: added
On the other hand, with $R_\alpha(\omega,u) =
\sup\{|\alpha- \beta|;(\beta,v)\stackrel{\omega}{\leftrightarrow
}(\alpha,u)\}$,
it is clear that (\ref{finite_2nd_moment}) holds when
$R_0(\omega,u)\leq\mbox{\textit{const}}$ for all
$u\in\omega_0$, ${\mathbb P}$-a.s. Such an example is given by the tube
$\{(\alpha,u)\in{\mathbb R}^2 \dvtx \cos\alpha\leq u \leq\cos\alpha
+1\}$,
a random shift to make
it stationary and ergodic (but not mixing).
\begin{rem}
(i) The continuity assumption of the map
$\alpha\mapsto\omega_\alpha$ has a
geometric meaning: it prevents the tube from having ``vertical
walls'' of nonzero surface measure. The reader may wonder what
happens without it. First, the disintegration formula
(\ref{differentials}) of the surface measure $\nu^\omega$ on
$\partial\omega$ becomes a product
$d {\bar\mu}^\omega_\alpha(u)\, d\phi^\omega(\alpha)$
where ${\bar\mu}^\omega_\alpha$ is a measure on the section of
$\partial\omega$
by the vertical hyperplane at $\alpha$ and where
$d \phi^\omega(\alpha)=\kappa_{\alpha,u}^{-1} \, d\alpha+
d \phi^\omega_{\mathrm{s}}
(\alpha)$ with a singular part $\phi^\omega_{\mathrm{s}}$.
If the singular part has atoms, one can see that the invariant
law ${\mathbb Q}$ [see (\ref{def_Q}) above] of the environment
seen from the particle has a marginal in $\omega$ which is singular
with respect to ${\mathbb P}$. This happens because, if the vertical walls
constitute a positive proportion of the tube's surface, in the
equilibrium the particle finds itself on a vertical wall with
positive probability;
on the other hand, if $\omega$ has the law ${\mathbb P}$,
a.s. there is no vertical wall at the origin.
The general situation is interesting but
complicated; in any case, our results continue to be valid in this
situation as well
[an important observation is that (\ref{IP<IQ})
below would still hold, possibly with another
constant]. To keep things simple, we will consider only,
all through the paper, random tubes satisfying the continuity
assumption. In the \hyperref[app]{Appendix}, we discuss the general case in more
detail. Another possible approach to this general case is to work
with the continuous-time stochastic billiard directly (cf.
Section \ref{s_env_continuous}).

(ii) A particular example of tubes is given by rotation invariant
tubes. They are obtained by rotating around the first axis the graph
of a positive bounded function. The main simplification is that,
with the proper formalism, one can forget the transverse
component $u$. Then the problem becomes purely one-dimensional.
\end{rem}

%s3 ###
\section{Environment viewed from the particle and the construction
of the corrector}
\label{s_environment}

%s3.1 ###
\subsection{Environment viewed from the particle: Discrete case}
\label{s_env_discrete}
One of the main methods we use in this paper is considering
the environment $\omega$ seen from the current location
of the random walk.
The ``environment viewed from the particle''
is the Markov chain
\[
\bigl((\theta_{\xi_n\cdot e}\omega,{\mathcal U}\xi_n), n=0,1,2,\ldots\bigr)
\]
with state space ${\mathfrak S}$.
The transition operator $G$ for this process
acts on functions $f\dvtx {\mathfrak S}\mapsto{\mathbb R}$ as follows
[cf. (\ref{differentials}) and (\ref{K->Phi})]:
%e13 ###
%
\begin{eqnarray}\label{def_trans_operator}
Gf (\omega,u) &=& {\mathtt E}_{\omega}\bigl(f(\theta_{\xi_1\cdot e}\omega
,{\mathcal U}\xi_1)
\mid\xi_0 = (0,u)\bigr)\nonumber\\
&=& \int_{\partial\omega} K ((0,u),x )
f(\theta_{x\cdot e}\omega,{\mathcal U}x) \, d\nu^{\omega}(x)
\\
&=& \int_{-\infty}^{+\infty} d\alpha\int_{\Lambda}
d\mu^{\omega}_\alpha(v)\, \kappa_{0,u}
f(\theta_\alpha\omega,v) \Phi(0,u,\alpha,v).\nonumber
\end{eqnarray}
\begin{rem}
Note that our environment consists not only of the tube with an
appropriate horizontal shift, but also of the transverse component
of the walk.
Another possible choice for the environment
would be obtained by rotating the shifted tube
to make it pass through the origin with inner normal at this point
given by the last coordinate vector. However, we made the present
choice to keep notation simple.
\end{rem}

Next, we show that this new Markov chain
is reversible with reversible
measure ${\mathbb Q}$ given by (\ref{def_Q}),
which means that $G$ is a self-adjoint operator in
$L^2({\mathfrak S})=L^2({\mathfrak S},{\mathbb Q})$:
\begin{lmm}
\label{l_revers_discr}
For all $f,g\in L^2({\mathfrak S})$ we have $\langle g, Gf \rangle
_{{\mathbb Q}} =
\langle f, Gg \rangle_{{\mathbb Q}}$.
Hence, the law ${\mathbb Q}$ is invariant for the Markov chain of the
environment viewed from the particle which means that for any
$f\in L^2({\mathfrak S})$ and all $n$,
%e14 ###
%
\begin{equation}
\label{IQ_is_invariant}
\langle{\mathtt E}_{\omega}[f(\theta_{\xi_n\cdot e}\omega
,{\mathcal U}\xi_n)\mid\xi
_0=(0,u)] \rangle_{{\mathbb Q}}
=\langle f \rangle_{{\mathbb Q}}.
\end{equation}
\end{lmm}
\begin{pf}
Indeed, from (\ref{def_Q}) and (\ref{def_trans_operator}),
%e18 ###
%e17 ###
%e16 ###
%e15 ###
%
\begin{eqnarray}\quad
\langle g, Gf \rangle_{{\mathbb Q}} &=& \frac{1}{{\mathcal Z}}\int
_{\Omega}d{\mathbb P}
\int_{\Lambda}
d\mu^{\omega}_0(u)\, g(\omega,u)\kappa_{0,u}^{-1}
\int_{-\infty}^{+\infty} d\alpha\nonumber\\
&&{} \times\int_{\Lambda} d\mu^{\omega}_\alpha(v)\, \kappa_{0,u}
\Phi(0,u,\alpha,v) f(\theta_\alpha\omega,v)\nonumber\\
\label{revers1}
&=& \frac{1}{{\mathcal Z}}\int_{-\infty}^{+\infty} d\alpha
\int_{\Omega}d{\mathbb P}
\int_{\Lambda^2} d\mu^{\omega}_{0,\alpha}(u,v)\, g(\omega,u)
f(\theta_\alpha\omega,v)\Phi(0,u,\alpha,v)
\\
\label{revers1'}
&=& \frac{1}{{\mathcal Z}}\int_{-\infty}^{+\infty} d\alpha
\int_{\Omega}d{\mathbb P}
\int_{\Lambda^2}
d\mu^{\theta_\alpha\omega}_{-\alpha,0}(u,v)\, g(\omega,u)
f(\theta_\alpha
\omega,v)\nonumber\\[-8pt]\\[-8pt]
&&\hspace*{96.7pt}{}\times\Phi^{\theta_\alpha\omega}(-\alpha,u,0,v)\nonumber\\
\label{revers2}
&=& \frac{1}{{\mathcal Z}}\int_{-\infty}^{+\infty} d\alpha
\int_{\Omega}d{\mathbb P}
\int_{\Lambda^2} d\mu^{\omega}_{-\alpha,0}(u,v)\,
g(\theta_{-\alpha} \omega,u)f(\omega,v)\nonumber\\[-8pt]\\[-8pt]
&&\hspace*{96.7pt}{}\times\Phi(-\alpha,u,0,v)\nonumber\\
&=& \frac{1}{{\mathcal Z}}\int_{\Omega}d{\mathbb P}
\int_{\Lambda} d\mu^{\omega}_0(v)\, f(\omega,v)\kappa_{0,v}^{-1}
\int_{-\infty}^{+\infty} d\alpha\nonumber\\
\label{revers3}
&&{} \times\int_{\Lambda} d\mu^{\omega}_{-\alpha}(u)\, \kappa_{0,v}
\Phi(0,v,-\alpha,u) g(\theta_{-\alpha} \omega,u)\\
&=& \langle f, Gg \rangle_{{\mathbb Q}},\nonumber
\end{eqnarray}
where we used (\ref{shift_Phi}) to pass from (\ref{revers1})
to (\ref{revers1'}),
the translation invariance of ${\mathbb P}$ to pass from (\ref{revers1'})
to (\ref{revers2}), the symmetry property (\ref{symmetry_Phi})
to pass from (\ref{revers2}) to (\ref{revers3}) and the
change of variable
$\alpha\mapsto- \alpha$ to obtain the last line.
\end{pf}

Let us define a semi-definite scalar product
$\langle g, f \rangle_{1}:=\langle g, (I-G)f \rangle_{{\mathbb Q}}$.
Again using (\ref{revers1}),
the translation invariance of ${\mathbb P}$ and the symmetry of $\Phi$,
we obtain
\begin{eqnarray*}
\langle g, f \rangle_{1}
&=& \frac{1}{{\mathcal Z}}\int_{-\infty}^{+\infty} d\alpha\int
_{\Omega}
d{\mathbb P}\int_{\Lambda^2} d\mu^{\omega}_{0,\alpha}(u,v)\,\Phi(0,u,\alpha,v)\\
&&\hspace*{97.74pt}{}\times
g(\omega,u) \bigl(f(\omega,u)-f(\theta_\alpha\omega,v)\bigr)\\
&=& \frac{1}{{\mathcal Z}}\int_{-\infty}^{+\infty} d\alpha\int
_{\Omega}
d{\mathbb P}
\int_{\Lambda^2} d\mu^{\omega}_{-\alpha,0}(u,v)\,
\Phi(-\alpha,u,0,v)\\
&&\hspace*{97.74pt}{}\times
g(\theta_{-\alpha}\omega,u)
\bigl(f(\theta_{-\alpha}\omega,u)-f(\omega,v)\bigr)\\
&=& \frac{1}{{\mathcal Z}}\int_{-\infty}^{+\infty} d\alpha
\int_{\Omega}d{\mathbb P}
\int_{\Lambda^2} d\mu^{\omega}_{0,\alpha}(u,v)\,
\Phi(0,u,\alpha,v)\\
&&\hspace*{97.74pt}{}\times
g(\theta_\alpha\omega,v)
\bigl(f(\theta_\alpha\omega,v)-f(\omega,u)\bigr).
\end{eqnarray*}
Consequently,
\begin{eqnarray*}
\langle g, f \rangle_{1} &=& \frac{1}{2{\mathcal Z}}\int_{-\infty
}^{+\infty}\hspace*{-1pt}
d\alpha\int_{\Omega}d{\mathbb P}
\int_{\Lambda^2} d\mu^{\omega}_{0,\alpha}(u,v)\,
\Phi(0,u,\alpha,v)
\bigl(f(\omega,u)-f(\theta_\alpha\omega,v)\bigr)\hspace*{1.5pt}\\
&&\hspace*{101.5pt}{} \times
\bigl(g(\omega,u)-g(\theta_\alpha\omega,v)\bigr),\hspace*{1.5pt}
\end{eqnarray*}
so the Dirichlet form $\langle f, f \rangle_{1}$ can be explicitly
written as
%e19 ###
%
\begin{eqnarray}
\label{eq_Dirichlet}
\langle f, f \rangle_{1} &=& \frac{1}{2{\mathcal Z}}\int_{-\infty
}^{+\infty}
d\alpha\int_{\Omega}d{\mathbb P}
\int_{\Lambda^2} d\mu^{\omega}_{0,\alpha}(u,v)\,
\Phi(0,u,\alpha,v)\nonumber\\[-8pt]\\[-8pt]
&&\hspace*{102.5pt}{}\times\bigl(f(\omega,u)-f(\theta_\alpha\omega,v)\bigr)^2,\nonumber
\end{eqnarray}
or, by (\ref{differentials}) and (\ref{K->Phi}),
%e20 ###
%
\begin{eqnarray}\label{eq_Dirichlet_K}
\langle f, f \rangle_{1} &=& \frac{1}{2{\mathcal Z}}\int_{\Omega
}d{\mathbb P}
\int_\Lambda\kappa_{0,u}^{-1} \,d\mu^{\omega}_0(u)\nonumber\\[-8pt]\\[-8pt]
&&{}\times\int_{\partial\omega} d\nu^{\omega}(x)\, K ((0,u),x )
\bigl(f(\omega,u)-f(\theta_{x\cdot e}\omega,{\mathcal
U}x)\bigr)^2.\nonumber
\end{eqnarray}

At this point it is convenient to establish the following result:
\begin{lmm}
\label{l_Q_ergodic}
The Markov process with initial law ${\mathbb Q}$ and transition
operator $G$ is ergodic.
\end{lmm}
\begin{pf}
Suppose that $f\in L^2({\mathfrak S})$
is such that $f=Gf$. Then $\langle f, f \rangle_{1}=0$
and so, by the translation invariance and (\ref{eq_Dirichlet_K}),
\[
\int_{\Omega}d{\mathbb P}\int_\Lambda\kappa_{s,u}^{-1} \,d\mu
^{\omega}_s(u)
\int_{\partial\omega} d\nu^{\omega}(x)\, K((s,u),x)
\bigl(f(\theta_s
\omega,u)-f(\theta_{x\cdot e}\omega,{\mathcal U}x)\bigr)^2 = 0
\]
for any $s$. Integrating the above equation in $s$ and
using (\ref{differentials}), we obtain
%e21 ###
%
\begin{equation}
\label{int_f2=0}
\int_{\Omega}d{\mathbb P}
\int_{(\partial\omega)^2} d\nu^{\omega}(x)\, d\nu^{\omega}(y)\, K(x,y)
\bigl(f(\theta_{x\cdot e}\omega,{\mathcal U}x)
-f(\theta_{y\cdot e}\omega,{\mathcal U}y) \bigr)^2 = 0.\hspace*{-30pt}
\end{equation}
Let us recall Lemma 3.3(iii) from \cite{CPSV1}:
if for some $x,y\in{\mathcal R}_{\omega}$ we have $K(x,y)>0$,
then there exist $\delta>0$ and two neighborhoods $B_x$ of $x$
and $B_y$ of $y$ such that
$K(x',y')>\delta$ for all $x'\in B_x, y'\in B_y$.
Now, for such $x,y$, (\ref{int_f2=0})
implies that there exists a constant ${\hat c}(\omega,x,y)$
such that $f(\theta_{z\cdot e}\omega,{\mathcal U}z)={\hat c}(\omega,x,y)$
for $\nu^{\omega}$-almost all $z\in B_x\cup B_y$. By the irreducibility
Condition \ref{ConditionP} (in fact, a much weaker
irreducibility condition would be sufficient), this implies that
$f(\theta_{z\cdot e}\omega,{\mathcal U}z)={\hat c}(\omega)$ for
$\nu^{\omega}$-almost all $z\in{\mathcal R}_{\omega}$.
Since ${\mathbb P}$ is ergodic, this means that $f(\omega,u)={\hat c}$
for $\mu^{\omega}_0$-almost all $u$ and ${\mathbb P}$-almost all
$\omega$.
\end{pf}

%s3.2 ###
\subsection{Environment viewed from the particle: Continuous case}
\label{s_env_continuous}
For the sake of completeness, we present also an analogous result
for the Knudsen stochastic billiard $(X_t,V_t)$. The notation and
the results of this section will not be used in the rest of this
paper.

Let $\widehat{\mathfrak S}= \{ (\omega,u)\dvtx \omega\in\Omega, u \in
\omega_0\}$,
and let ${\widehat{\mathbb Q}}$ be the probability measure
on $\widehat{\mathfrak S}\times{\mathbb S}^{d-1}$
defined by
\[
d{\widehat{\mathbb Q}}(\omega,u,h) = \frac{1}{{\widehat{\mathcal Z}}}
{\mathbb I}{\{u\in\omega_0\}} \,du \,dh \,d{\mathbb P}(\omega),
\]
where ${\widehat{\mathcal Z}}=|{\mathbb S}^{d-1}|\int_\Omega|\omega
_0|d{\mathbb P}$
is the normalizing constant.
The scalar product $\langle\cdot, \cdot\rangle_{{\widehat{\mathbb
Q}}}$ in
$L^2({\widehat{\mathbb Q}})$ is given,
for two ${\widehat{\mathbb Q}}$-square integrable functions
${\hat f}, {\hat g}\dvtx \widehat{\mathfrak S}\times{\mathbb S}^{d-1}
\mapsto{\mathbb R}$,
by
\[
\langle{\hat f}, {\hat g} \rangle_{{\widehat{\mathbb Q}}} = \frac
{1}{{\widehat{\mathcal Z}}}
\int_\Omega d{\mathbb P}\int_{\omega_0}du
\int_{{\mathbb S}^{d-1}}dh\,
{\hat f}(\omega,u,h){\hat g}(\omega,u,h).
\]
For the continuous\vspace*{1pt} time KSB, the ``environment viewed
from the particle'' is the Markov process
$ ((\theta_{X_t\cdot e}\omega,{\mathcal U}X_t,V_t),t\geq0 )$
with the state space
$\widehat{\mathfrak S}\times{\mathbb S}^{d-1}$.
The transition semi-group ${\widehat{\mathcal P}}_t$ for
this process acts on functions
${\hat f}\dvtx \widehat{\mathfrak S}\times{\mathbb S}^{d-1} \mapsto
{\mathbb R}$ in the following
way:
\[
{\widehat{\mathcal P}}_t{\hat f}(\omega,u,h)
= {\mathtt E}_{\omega}\bigl({\hat f}(\theta_{X_t\cdot e}\omega,{\mathcal
U}X_t,V_t)
\mid X_0=(0,u), V_0=h \bigr).
\]
We show that ${\widehat{\mathcal P}}$ is
quasi-reversible with respect to the law ${\widehat{\mathbb Q}}$.
\begin{lmm}
\label{l_revers_cont}
For all ${\hat f}, {\hat g}\in L^2({\widehat{\mathbb Q}})$ and $t>0$
we have
%e22 ###
%
\begin{equation}
\label{eq:quasirev}
\langle{\hat g}, {\widehat{\mathcal P}}_t{\hat f} \rangle_{{\widehat
{\mathbb Q}}} =
\langle{\breve f}, {\widehat{\mathcal P}}_t{\breve g} \rangle
_{{\widehat{\mathbb Q}}}
\end{equation}
with $\breve f (\omega,u,h)=\hat f (\omega,u,-h)$.
In particular, the law ${\widehat{\mathbb Q}}$ is invariant for the process
$ ((\theta_{X_t\cdot e}\omega,{\mathcal U}X_t,V_t),t\geq0 )$.
\end{lmm}
\begin{pf}
We first prove that (\ref{eq:quasirev}) implies that the
law ${\widehat{\mathbb Q}}$ is invariant. Indeed, taking ${\hat
g}={\mathbf1}$,
we get for all test functions ${\hat f}$
\begin{eqnarray*}
\langle{\mathbf1}, {\widehat{\mathcal P}}_t{\hat f} \rangle
_{{\widehat{\mathbb Q}}} &=&
\langle{\breve f}, {\mathbf1} \rangle_{{\widehat{\mathbb Q}}}\\
&=& \langle{ f}, {\mathbf1} \rangle_{{\widehat{\mathbb Q}}}
\end{eqnarray*}
by the change of variable $h$ into $-h$ in the integral.
Hence ${\widehat{\mathbb Q}}$ is invariant.

We now turn to the proof of (\ref{eq:quasirev}).
Introducing the notation ${\mathtt P}_t^{\omega}$
for the transition kernel of KSB,
\[
{\mathtt P}_{\omega}( X_t \in dx', V_t \in dh' \mid X_0=x, V_0=h) =
{\mathtt P}_t^{\omega}( x,h; x',h') \,dx' \,dh',
\]
we observe that
\[
{\widehat{\mathcal P}}_t{\hat f}(\omega,u,h) =
\int_{\omega}dx'
\int_{{\mathbb S}^{d-1}}dh'\,
{\mathtt P}_t^{\omega}( (0,u),h;x',h' )
{\hat f}(\theta_{x' \cdot e}, \omega, {\mathcal U}x',h').
\]
In Theorem 2.5 in \cite{CPSV1}, it was shown that ${\mathtt
P}_t^{\omega}$
is itself quasi-reversible, that is,
\[
{\mathtt P}_t^{\omega}(x,h;x',h')={\mathtt P}_t^{\omega}(x',-h';x,-h).
\]
Therefore,
\begin{eqnarray*}
\langle{\hat g}, {\widehat{\mathcal P}}_t{\hat f} \rangle_{{\widehat
{\mathbb Q}}}
&=&
\frac{1}{{\widehat{\mathcal Z}}} \int_\Omega d{\mathbb P}\int
_{\omega_0}du
\int_{{\mathbb S}^{d-1}}dh\, {\hat g}(\omega,u,h)
{\widehat{\mathcal P}}_t {\hat f}(\omega,u,h)
\\
&=&
\frac{1}{{\widehat{\mathcal Z}}} \int_\Omega d{\mathbb P}\int
_{\omega_0}du
\int_{{\mathbb S}^{d-1}}dh\, {\hat g}(\omega,u,h)
\int_{\mathbb R}d \alpha
\int_{\omega_\alpha}du' \\
&&{} \times\int_{{\mathbb S}^{d-1}}dh'\,
{\mathtt P}_t^{\omega}( (0,u),h;(\alpha,u'),h' )
{\hat f}(\theta_{\alpha}, \omega, u',h')
\\
&=&
\frac{1}{{\widehat{\mathcal Z}}} \int_\Omega d{\mathbb P}\int
_{\omega_0}du
\int_{{\mathbb S}^{d-1}}dh\, {\hat g}(\omega,u,h)
\int_{\mathbb R}d \alpha
\int_{\omega_\alpha}du' \\
&&{} \times\int_{{\mathbb S}^{d-1}}dh'\,
{{\mathtt P}_t^{\theta_\alpha\omega}{}}
( (0,u'),-h';(-\alpha,u),-h )
{\hat f}(\theta_{\alpha} \omega, u',h')
\\
&=&
\frac{1}{{\widehat{\mathcal Z}}}
\int_{\mathbb R}d \alpha
\int_\Omega d{\mathbb P}\int_{\omega_{-\alpha}}du
\int_{{\mathbb S}^{d-1}}dh\, {\hat g}(\theta_{-\alpha},\omega,u,h)
\int_{\omega_0}du' \\
&&{} \times\int_{{\mathbb S}^{d-1}}dh'\,
{{\mathtt P}_t^{\omega}{}}
( (0,u'),-h';(-\alpha,u),-h )
{\hat f}( \omega, u',h')
\\
&=&
\frac{1}{{\widehat{\mathcal Z}}}
\int_{\mathbb R}d \alpha
\int_\Omega d{\mathbb P}\int_{\omega_{\alpha}}du
\int_{{\mathbb S}^{d-1}}dh\, {\hat g}(\theta_{\alpha},\omega,u,-h)
\int_{\omega_0}du' \\
&&{} \times
\int_{{\mathbb S}^{d-1}}dh'\,
{{\mathtt P}_t^{\omega}{}}
( (0,u'),h';(\alpha,u),h )
{\hat f}( \omega, u',-h')
\\&=& \langle{\breve f}, {\widehat{\mathcal P}}_t{\breve g} \rangle
_{{\widehat{\mathbb Q}}},
\end{eqnarray*}
where we used that the Lebesgue measure on ${\mathbb R}^d$ is product
to get the second line, quasi-reversibility for the third one,
Fubini and translation invariance of ${\mathbb P}$ for the fourth one, and
change of variables $(h,h',\alpha)$ to $(-h,-h',-\alpha)$ in the
fifth one.
%\rightqed
\end{pf}

%s3.3 ###
\subsection{Construction of the corrector}
\label{s_constr_corrector}
Now, we are going to construct the corrector function for the
random walk $\xi$.

Let us show that for any $g\in L^2({\mathfrak S})$,
%e23 ###
%
\begin{equation}
\label{test_function}
\langle g, \Delta\rangle_{{\mathbb Q}} \leq\tfrac{1}{\sqrt
{2}}\langle b
\rangle_{{\mathbb Q}}^{1/2}\langle g, g \rangle_{1}^{1/2}.
\end{equation}
Indeed, from (\ref{def_drift}) we obtain
\begin{eqnarray*}
\langle g, \Delta\rangle_{{\mathbb Q}} &=& \frac{1}{{\mathcal Z}}\int
_{-\infty
}^{+\infty} \alpha\,
d\alpha\int_\Omega d{\mathbb P}
\int_{\Lambda^2}\mu^{\omega}_{0,\alpha}(u,v) g(\omega,u)
\Phi(0,u,\alpha,v)\\
&=& \frac{1}{{\mathcal Z}}\int_{-\infty}^{+\infty} \alpha
\,d\alpha\int_\Omega d{\mathbb P}
\int_{\Lambda^2}\mu^{\omega}_{-\alpha,0}(u,v)
g(\theta_{-\alpha}\omega,u)\Phi(-\alpha,u,0,v)\\
&=& -\frac{1}{{\mathcal Z}}\int_{-\infty}^{+\infty} \alpha
\,d\alpha\int_\Omega d{\mathbb P}
\int_{\Lambda^2}\mu^{\omega}_{0,\alpha}(u,v)
g(\theta_\alpha\omega,v)\Phi(0,u,\alpha,v),
\end{eqnarray*}
so
%e24 ###
%
\begin{eqnarray}
\label{g_Delta}
\langle g, \Delta\rangle_{{\mathbb Q}} &=& \frac{1}{2{\mathcal Z}}\int
_{-\infty
}^{+\infty}
\alpha \,d\alpha\int_\Omega d{\mathbb P}\nonumber\\[-8pt]\\[-8pt]
&&{}\times\int_{\Lambda^2}\mu^{\omega}_{0,\alpha}(u,v)
\Phi(0,u,\alpha,v)\bigl(g(\omega,u)-g(\theta_\alpha\omega,v)\bigr).\nonumber
\end{eqnarray}
Using the Cauchy--Schwarz inequality in (\ref{g_Delta}), we obtain
\begin{eqnarray*}
\langle g, \Delta\rangle_{{\mathbb Q}} &\leq& \frac{1}{2}
\biggl[ \frac{1}{{\mathcal Z}}\int_{-\infty}^{+\infty} \alpha^2
\,d\alpha\int_\Omega d{\mathbb P}
\int_{\Lambda^2}\mu^{\omega}_{0,\alpha}(u,v) \Phi(0,u,\alpha,v)\\
&&\hspace*{11pt}{} \times\frac{1}{{\mathcal Z}}\int_{-\infty}^{+\infty} d\alpha
\int_\Omega d{\mathbb P}
\int_{\Lambda^2}\mu^{\omega}_{0,\alpha}(u,v) \Phi(0,u,\alpha,v)
\\
&&\hspace*{120.6pt}{}\times\bigl(g(\omega,u)-g(\theta_\alpha\omega,v)\bigr)^2 \biggr]^{1/2}\\
&=& \frac{1}{2} \langle b \rangle_{{\mathbb Q}}^{1/2}(2\langle g, g
\rangle
_{1})^{1/2},
\end{eqnarray*}
which shows (\ref{test_function}).

Note that, from (\ref{g_Delta}) with $g=1$ we obtain the
important property
\[
\langle\Delta\rangle_{{\mathbb Q}}=0.
\]
As shown in Chapter 1 of \cite{KLO}, we have the variational
formula
\begin{eqnarray*}
\langle(I-G)^{-1/2}\Delta, (I-G)^{-1/2}\Delta\rangle_{{\mathbb Q}} &=&
\langle\Delta, (I-G)^{-1}\Delta\rangle_{{\mathbb Q}}\\
&=&\sup\{\langle g, \Delta\rangle_{{\mathbb Q}} - \langle g, g
\rangle
_{1}, \langle g, g \rangle_{1}<\infty\}.
\end{eqnarray*}
Then provided that (\ref{finite_2nd_moment}) holds, inequality (\ref
{test_function}) implies that the
drift $\Delta$ belongs to the range of
$(I-G)^{1/2}$, and so the time-variance of $\Delta$ is finite.
%FC: At this point we mention that this
% already implies the annealed CLT,
% see \cite{KV}, \cite{DFGW}.
At this point we mention that this already implies weaker
forms of the CLT, by applying \cite{KV}
(under the invariant measure,
or in probability with respect to the environment)
or \cite{DFGW} (under the annealed measure).
With this observation, we could have used the resolvent
method originally developed in \cite{KV,KLO}
to construct the corrector.
However, it is more direct to use the method of the orthogonal
projections on the potential subspace (cf. \cite{BP,M,MP}).

For $\omega\in\Omega$, $u\in\partial\omega_0$, define
\begin{eqnarray*}
{\mathsf V}^+_{\omega,u} &=& \{y\in\partial\omega\dvtx y\cdot e > 0,
K ((0,u),y )>0\},\\
{\mathsf V}^-_{\omega,u} &=& \{y\in\partial\omega\dvtx y\cdot e < 0,
K ((0,u),y )>0\}.
\end{eqnarray*}
Then, in addition to the space ${\mathfrak S}$, we define two spaces
${\mathfrak N}\subset{\mathfrak M}$ in the following way:
\begin{eqnarray*}
{\mathfrak N}&=& \{(\omega,u,y) \dvtx \omega\in\Omega,u\in\partial
\omega_0,
y\in{\mathsf V}^+_{\omega,u}\},\\
{\mathfrak M}&=& \{(\omega,u,y) \dvtx \omega\in\Omega,u\in\partial
\omega_0,
y\in\partial\omega\}.
\end{eqnarray*}
On ${\mathfrak N}$ we define the measure $K {\mathbb Q}$ with mass that
is less than 1
for which
a nonnegative function $f\dvtx{\mathfrak N}\mapsto{\mathbb R}$ has the
expected value
%e25 ###
%
\begin{equation}
\label{df_ex_KQ}
\langle f \rangle_{{K{\mathbb Q}}} = \biggl\langle\int_{{\mathsf
V}^+_{\omega
,u}}f(\omega,u,y) K ((0,u),y )\, d\nu^{\omega}(y) \biggr\rangle_{{\mathbb Q}}.
\end{equation}
Two square-integrable functions $f,g\dvtx{\mathfrak N}\mapsto{\mathbb R}$ have
scalar product,
%e26 ###
%
\begin{equation}
\label{df_sc_KQ}
\langle f, g \rangle_{{K{\mathbb Q}}} = \biggl\langle\int_{{\mathsf
V}^+_{\omega
,u}}f(\omega,u,y) g(\omega,u,y) K ((0,u),y ) \,d\nu^{\omega}(y)
\biggr\rangle_{{\mathbb Q}}.
\end{equation}

Also, define the gradient $\nabla$ as the map that transfers a
function $f\dvtx {\mathfrak S}\mapsto{\mathbb R}$
to the function $\nabla f \dvtx {\mathfrak N}\mapsto{\mathbb R}$, defined by
%e27 ###
%
\begin{equation}
\label{df_nabla}
(\nabla f)(\omega,u,y) = f(\theta_{y\cdot e}\omega, {\mathcal U}y)
- f(\omega,u).
\end{equation}
Since ${\mathbb Q}$ is reversible for $G$, we can write
\begin{eqnarray*}
\langle(\nabla f)^2 \rangle_{{K{\mathbb Q}}} &=& \biggl\langle\int
_{{\mathsf V}
^+_{\omega,u}} \bigl(f(\theta_{y\cdot e}\omega, {\mathcal U}y) - f(\omega
,u) \bigr)^2
K ((0,u),y )\, d\nu^{\omega}(y) \biggr\rangle_{{\mathbb Q}}\\
&\leq& 2\biggl\langle\int_{\partial\omega} f^2(\theta_{y\cdot e}\omega
, {\mathcal U}y) K ((0,u),y ) \,d\nu^{\omega}(y) \biggr\rangle_{{\mathbb Q}}
\\
&&{} + 2\biggl\langle\int_{\partial\omega}f^2(\omega, u) K ((0,u),y )
\,d\nu^{\omega}(y) \biggr\rangle_{{\mathbb Q}}\\
&=& 2\langle Gf^2 \rangle_{{\mathbb Q}} + 2\langle f^2 \rangle
_{{\mathbb Q}}\\
&=& 4\langle f^2 \rangle_{{\mathbb Q}},
\end{eqnarray*}
so $\nabla$ is, in fact, a map from $L^2({\mathfrak S})$ to
$L^2({\mathfrak N})$.

Then, following \cite{BP}, we denote by $L^2_\nabla({\mathfrak N})$
the closure of the set of gradients of all
functions from $L^2({\mathfrak S})$.
We then consider the orthogonal decomposition of
$L^2({\mathfrak N})$
into the
``potential'' and the ``solenoidal'' subspaces:
$L^2({\mathfrak N})=L^2_\nabla({\mathfrak N})\oplus(L^2_\nabla
({\mathfrak N}))^\bot$. To
characterize the solenoidal subspace $(L^2_\nabla({\mathfrak N}))^\bot
$, we
introduce the divergence operator
in the following way. For $f\dvtx{\mathfrak N}\mapsto{\mathbb R}$, we have
$\operatorname{div}f \dvtx {\mathfrak S}
\mapsto{\mathbb R}$
defined by
%e28 ###
%
\begin{eqnarray}\label{df_divergence}
(\operatorname{div}f)(\omega,u) &=& \int_{{\mathsf V}^+_{\omega,u}}K
((0,u),y )
f(\omega,u,y) \,d\nu^{\omega}(y)\nonumber\\[-8pt]\\[-8pt]
&&{} - \int_{{\mathsf V}^-_{\omega,u}}K ((0,u),y )
f \bigl(\theta_{y\cdot e}\omega,{\mathcal U}y,(|y\cdot e|,u) \bigr)\,
d\nu^{\omega}(y)\nonumber
\end{eqnarray}
[note that for $y\in{\mathsf V}^-_{\omega,u}$ we have
$(|y\cdot e|,u)\in{\mathsf V}^+_{\theta_{y\cdot e}\omega,{\mathcal U}y}$].
Now, we verify the following integration by parts formula:
for any $f\in L^2({\mathfrak S})$, $g\in L^2({\mathfrak N})$,
%e29 ###
%
\begin{equation}
\label{gradient=-diverg}
\langle g, \nabla f \rangle_{{K{\mathbb Q}}} = -\langle f\operatorname
{div}g \rangle_{{\mathbb Q}}.
\end{equation}
Indeed, we have
%e30 ###
%
\begin{eqnarray}\label{calc_g_nabla_f}\quad
\langle g, \nabla f \rangle_{{K{\mathbb Q}}} &=& \biggl\langle\int
_{{\mathsf V}
^+_{\omega,u}} K ((0,u),y )g(\omega,u,y)\nonumber\\
&&\hspace*{24.3pt}{}\times \bigl(f(\theta_{y\cdot e}\omega
,{\mathcal U}y) -f(\omega,u) \bigr) \,d\nu^{\omega}(y) \biggr\rangle_{{\mathbb
Q}}\nonumber\\[-8pt]\\[-8pt]
&=& - \biggl\langle\int_{{\mathsf V}^+_{\omega,u}} K ((0,u),y ) g(\omega
,u,y)f(\omega,u)\, d\nu^{\omega}(y) \biggr\rangle_{{\mathbb Q}}
\nonumber\\
&&{} + \biggl\langle\int_{{\mathsf V}^+_{\omega,u}} K ((0,u),y ) g(\omega,u,y)
f(\theta_{y\cdot e}\omega,{\mathcal U}y) \, d\nu^{\omega}(y) \biggr\rangle
_{{\mathbb Q}}.\nonumber
\end{eqnarray}
For the second term in the right-hand side
of (\ref{calc_g_nabla_f}), we obtain
\begin{eqnarray*}
&&\biggl\langle\int_{{\mathsf V}^+_{\omega,u}} K ((0,u),y )
g(\omega
,u,y) f(\theta_{y\cdot e}\omega,{\mathcal U}y) \,d\nu^{\omega}(y)
\biggr\rangle_{{\mathbb Q}}\\
&&\qquad= \frac{1}{{\mathcal Z}}\int_\Omega d{\mathbb P}\int_0^{+\infty
}d\alpha
\int_{\Lambda^2} d\mu^{\omega}_{0,\alpha}(u,v)\,
\Phi(0,u,\alpha,v)g (\omega,u,(\alpha,v) )
f(\theta_\alpha\omega,v)\\
&&\qquad= \frac{1}{{\mathcal Z}}\int_{-\infty}^0 d\alpha\int_\Omega
d{\mathbb P}
\int_{\Lambda^2} d\mu^{\omega}_{\alpha,0}(v,u)\,
\Phi(\alpha,v,0,u)g (\theta_\alpha\omega,v,(|\alpha|,u) )
f(\omega,u)\\
&&\qquad= \biggl\langle\int_{{\mathsf V}^-_{\omega,u}}f(\omega,u) g \bigl(\theta
_{y\cdot
e}\omega,{\mathcal U}y,(|y\cdot e|,u) \bigr) K ((0,u),y ) \,d\nu^{\omega
}(y) \biggr\rangle_{{\mathbb Q}},
\end{eqnarray*}
and so
\begin{eqnarray*}
\langle g, \nabla f \rangle_{{K{\mathbb Q}}} &=& - \biggl\langle\int
_{{\mathsf V}
^+_{\omega,u}} g(\omega,u,y)f(\omega,u)K ((0,u),y ) \,d\nu^{\omega}(y)
\biggr\rangle_{{\mathbb Q}}\\
&&{}+ \biggl\langle\int_{{\mathsf V}^-_{\omega,u}}f(\omega,u) g \bigl(\theta
_{y\cdot
e}\omega,{\mathcal U}y,(|y\cdot e|,u) \bigr) K ((0,u),y ) \,d\nu^{\omega
}(y) \biggr\rangle_{{\mathbb Q}
}\\
&=& -\langle f\operatorname{div}g \rangle_{{\mathbb Q}}
\end{eqnarray*}
and the proof of (\ref{gradient=-diverg}) is complete.

Analogously to Lemma 4.2 of \cite{BP}, we can characterize the
space $(L^2_\nabla({\mathfrak N}))^\bot$ as follows:
\begin{lmm}
\label{l_charact_solenoidal}
$g\in(L^2_\nabla({\mathfrak N}))^\bot$ if and only if $\operatorname
{div}g(\omega,u) = 0$
for ${\mathbb Q}$-almost all $(\omega,u)$.
\end{lmm}
\begin{pf}
This is a direct consequence of (\ref{gradient=-diverg}).
\end{pf}

A function $f\in L^2({\mathfrak N})$ can be interpreted as a flow by
putting formally
\[
f(\omega,u,y) := -
f \bigl(\theta_{y\cdot e}\omega,{\mathcal U}y, (|y\cdot e|,u) \bigr)
\]
for $y\in{\mathsf V}^-_{\omega,u}, \omega\in{\mathfrak S}$.
Then it is straightforward to obtain that every $f\in
L^2_\nabla({\mathfrak N})$ is \textit{curl-free}, which means that for
any loop
$y_0, y_1,\ldots,y_n\in\partial\omega$ with $y_0=y_n$ and
$K(y_i,y_{i+1})>0$ for $i=1,\ldots,n-1$, we have
%e31 ###
%
\begin{equation}
\label{cocycle}
\sum_{i=0}^{n-1} f \bigl(\theta_{y_i\cdot e}\omega,{\mathcal U}y_i,
y_{i+1}-(y_i\cdot e)e \bigr) = 0.
\end{equation}
Every curl-free function $f$ can be integrated into a unique
function $\phi\dvtx {\mathfrak M}\mapsto{\mathbb R}$
which can be defined by
%e32 ###
%
\begin{equation}
\label{assembling_phi}
\phi(\omega,u,y) = \sum_{i=0}^{n-1}
f \bigl(\theta_{y_i\cdot e}\omega,{\mathcal U}y_i,
y_{i+1}-(y_i\cdot e)e \bigr),
\end{equation}
where $y_0, y_1,\ldots,y_n \in\partial\omega$ is an arbitrary path
such that
$y_0=(0,u)$, $y_n=y$, and $K(y_i,y_{i+1})>0$ for $i=1,\ldots,n-1$.
Automatically, such a function $\phi$ satisfies the following
\textit{shift-covariance} property: for any $u\in\partial\omega_0$,
$x,y\in\partial\omega$,
%e33 ###
%
\begin{equation}
\label{shift-covariance}
\phi(\omega,u,x) = \phi(\omega,u,y)
+ \phi\bigl(\theta_{y\cdot e}\omega,{\mathcal U}y,x-(y\cdot e)e \bigr).
\end{equation}
We denote by ${\mathcal H}({\mathfrak M})$ the set of all
shift-covariant functions
${\mathfrak M}\to{\mathbb R}$.
Note that, by taking $x=y=(0,u)$ in (\ref{shift-covariance}),
we obtain
%e34 ###
%
\begin{equation}
\label{shift-cov_in_0}
\phi(\omega,u,(0,u) )=0 \qquad\mbox{for any } \phi\in
{\mathcal H}({\mathfrak M}).
\end{equation}
Also, for any $\phi\in{\mathcal H}({\mathfrak M})$, we define
$\operatorname{grad}\phi$ as the
unique function $f\dvtx{\mathfrak N}\to{\mathbb R}$,
from the shifts of which $\phi$ is assembled
[as in (\ref{assembling_phi})].
In view of (\ref{shift-cov_in_0}), we can write
\[
(\operatorname{grad}f)(\omega,u,y) = f(\omega,u,y) \qquad\mbox{for }
(\omega,u,y)\in{\mathfrak N}, f \in{\mathcal H}({\mathfrak M}).
\]

Let us define
an operator ${\mathcal L}$ which transfers a function $\phi\dvtx{\mathfrak
M}\mapsto{\mathbb R}$
to a function
$f\dvtx{\mathfrak S}\mapsto{\mathbb R}$, $f={\mathcal L}\phi$ with
%e35 ###
%
\begin{equation}
\label{df_action_G_N}\quad
({\mathcal L}\phi)(\omega,u) = \int_{\partial\omega}
K ((0,u),y ) [\phi(\omega,u,y)-\phi(\omega,u,(0,u))
] \, d\nu^{\omega}(y).
\end{equation}
Note that, by (\ref{df_nabla}), we obtain ${\mathcal L}(\nabla f) = Gf-f$
for any $f \in L^2({\mathfrak S})$.
Then, using (\ref{shift-covariance}) and (\ref{shift-cov_in_0}),
we write, for $\phi\in{\mathcal H}({\mathfrak M})$,
\begin{eqnarray*}
(\operatorname{div}\operatorname{grad}\phi)(\omega,u) &=& \int
_{{\mathsf V}^+_{\omega,u}}K ((0,u),y
)\phi(\omega,u,y) \,d\nu^{\omega}(y)\\
&&{}-\int_{{\mathsf V}^-_{\omega,u}}K ((0,u),y )
\phi\bigl(\theta_{y\cdot e}\omega,{\mathcal U}y,(|y\cdot e|,u)\bigr) \,d\nu^{\omega}(y)\\
&=& \int_{{\mathsf V}^+_{\omega,u}\cup{\mathsf V}^-_{\omega,u}}K ((0,u),y
)\phi(\omega,u,y) \,d\nu^{\omega}(y).
\end{eqnarray*}
So, for any $\phi\in{\mathcal H}({\mathfrak M})$, we have
$\operatorname{div}\operatorname{grad}\phi= {\mathcal L}\phi$.
This observation together with Lemma \ref{l_charact_solenoidal}
immediately implies the following fact:
\begin{lmm}
\label{l_harmonic}
Suppose that $\phi\in{\mathcal H}({\mathfrak M})$ is such that
$\operatorname{grad}\phi\in(L^2_\nabla({\mathfrak N}))^\bot$.
Then $\phi$ is harmonic for the Knudsen random walk, that is,
$({\mathcal L}\phi)(\omega,u)=0$ for ${\mathbb Q}$-almost all
$(\omega,u)$.
\end{lmm}

Now, we are able to construct the corrector.
Consider the function\break $\phi(\omega,u,y)=y\cdot e$,
and observe that $\phi\in{\mathcal H}({\mathfrak M})$.
Let ${\hat\phi}=\operatorname{grad}\phi$. First, let us show that
%e36 ###
%
\begin{equation}
\label{hatphi=b/2}
\langle{\hat\phi}^2 \rangle_{{K{\mathbb Q}}} = \tfrac{1}{2}\langle b
\rangle_{{\mathbb Q}},
\end{equation}
that is, if (\ref{finite_2nd_moment}) holds, then
${\hat\phi}\in L^2({\mathfrak N})$. Indeed,
%e37 ###
%
\begin{eqnarray}\label{calculation_symmetry}
\langle{\hat\phi}^2 \rangle_{{K{\mathbb Q}}} &=& \biggl\langle\int
_{{\mathsf V}
^+_{\omega,u}} (y\cdot e)^2 K ((0,u),y ) \,d\nu^{\omega}(y) \biggr\rangle
_{{\mathbb Q}
}\nonumber\\
&=& \frac{1}{{\mathcal Z}}\int_0^{+\infty}\alpha^2\,
d\alpha\int_\Omega d{\mathbb P}\int_{\Lambda^2}
d\mu^{\omega}_{0,\alpha}(u,v)\,
\Phi(0,u,\alpha,v)\nonumber\\[-8pt]\\[-8pt]
&=& \frac{1}{{\mathcal Z}}\int_{-\infty}^0\alpha^2\,
d\alpha\int_\Omega d{\mathbb P}\int_{\Lambda^2}d
\mu^{\omega}_{\alpha,0}(v,u)\,
\Phi(\alpha,v,0,u)\nonumber\\
&=& \biggl\langle\int_{{\mathsf V}^-_{\omega,u}}(y\cdot e)^2 K ((0,u),y )\,
d\nu^{\omega}
(y) \biggr\rangle_{{\mathbb Q}}\nonumber
\end{eqnarray}
and so
\begin{eqnarray*}
\langle{\hat\phi}^2 \rangle_{{K{\mathbb Q}}} &=& \frac{1}{2}\biggl\langle
\int
_{{\mathsf V}^+_{\omega,u} \cup{\mathsf V}^-_{\omega,u}}(y\cdot e)^2
K ((0,u),y )\,
d\nu^{\omega}(y) \biggr\rangle_{{\mathbb Q}}\\
&=& \frac{1}{2}\langle b \rangle_{{\mathbb Q}}.
\end{eqnarray*}
Then, let $g$ be the orthogonal projection of $(-{\hat\phi})$
onto $L^2_\nabla({\mathfrak N})$.
Define $\psi\in{\mathcal H}({\mathfrak M})$ to be the unique function
such that
$g=\operatorname{grad}\psi$;
in particular,\break $\psi(\omega,u,(0,u) )=0$ for
$u\in\partial\omega_0$.
Then we have
\[
{\hat\phi}+g = \operatorname{grad}\bigl((y\cdot e)+\psi(\omega,u,y) \bigr)
\in
(L^2_\nabla({\mathfrak N}))^\bot,
\]
so Lemma \ref{l_harmonic} implies that for ${\mathbb Q}$-a.a.
$(\omega,u)$, $\psi$ is the corrector in the sense that
%e38 ###
%
\begin{equation}
\label{psi_is_cr}
{\mathtt E}_{\omega}\bigl((\xi_1-\xi_0)\cdot e + \psi(\omega,u,\xi_1)
-\psi(\omega,u,\xi_0)\mid\xi_0=(0,u) \bigr)=0
\end{equation}
[recall that, by (\ref{shift-cov_in_0}), the term
$\psi(\omega,u,\xi_0)$ can be dropped from (\ref{psi_is_cr})].
Now, denote
\[
J^\omega_x = {\mathtt E}_{\omega}\bigl((\xi_1-\xi_0)\cdot e
+ \psi\bigl(\theta_{\xi_0\cdot e}\omega,{\mathcal U}\xi_0,
\xi_1-(\xi_0\cdot e)e\bigr)\mid\xi_0=x \bigr).
\]
By the translation invariance of ${\mathbb P}$, (\ref{psi_is_cr})
and (\ref{differentials}), we can write
\begin{eqnarray*}
0 &=& \int_{-\infty}^{+\infty}d\alpha
\frac{1}{{\mathcal Z}}\int_\Omega d{\mathbb P}\int_{\Lambda}d\mu
^{\omega}_\alpha(u)\,
\kappa^{-1}_{\alpha,u}\bigl|J^\omega_{(\alpha,u)}\bigr|\\
&=& \frac{1}{{\mathcal Z}}\int_\Omega d{\mathbb P}\int_{\partial
\omega}
|J^\omega_x| \,d\nu^{\omega}(x)
\end{eqnarray*}
and this implies that $J_x^\omega=0$ for
$\nu^\omega\otimes{\mathbb P}$-a.e.
%FC: $(\omega,u)$.
$(\omega,x)$.
From (\ref{shift-covariance}), we have
\[
\psi(\omega,u,y)-\psi(\omega,u,x)=
\psi\bigl(\theta_{x\cdot e}\omega,{\mathcal U}x,y-(x\cdot e)e\bigr),
\]
which does not depend on $u$. We summarize this in:
\begin{prop}
\label{Prop_corrector}
For ${\mathbb P}$-almost all $\omega$, it holds
%e39 ###
%
\begin{equation}
\label{psi_is_corrector}
{\mathtt E}_{\omega}\bigl((\xi_1-\xi_0)\cdot e + \psi(\omega,u,\xi_1)
-\psi(\omega,u,\xi_0)\mid\xi_0=x \bigr)=0
\end{equation}
for all $u\in\partial\omega_0$
and $\nu^{\omega}$-almost all $x\in\partial\omega$.
\end{prop}

%s3.4 ###
\subsection{Sequence of reference points and properties of
the corrector}
\label{s_refpoints}
Let $\chi$ be a random variable with uniform distribution in
$[-1/2,1/2]$,
independent of everything. Note that $(\chi+n, n\in{\mathbb Z})$
is then a stationary point process on the real line.
For a fixed environment $\omega$ define the sequence of
conditionally
independent random variables $\zeta_n\in\Lambda$, $n\in{\mathbb Z}$,
with $\zeta_n$ uniformly distributed on $\partial\omega_{\chi+n}$,
%e40 ###
%
\begin{equation}
\label{distr_zeta}
{\mathtt P}_{\omega}[\zeta_n\in B] = \frac{\mu^{\omega}_{\chi
+n}(B)}{\mu^{\omega}_{\chi+n}
(\partial\omega_{\chi+n})}.
\end{equation}
We denote by $\mathrm{E}^{\zeta}$ the expectation with respect to
$\zeta$
and $\chi$ (with fixed $\omega$),
and by $\bar{\mathrm{E}}^{\zeta}$ the expectation with respect to
$\zeta$
conditioned on $\{\chi=0\}$.
Then by (\ref{shift-covariance})
we have the following decomposition:
%e41 ###
%
\begin{equation}
\label{decomp_corr}
\psi(\theta_\chi\omega,\zeta_0,(n,\zeta_n) ) =
\sum_{i=0}^{n-1}\psi(\theta_{\chi+i}\omega,\zeta_i,
(1,\zeta_{i+1}) )
\end{equation}
so that $\psi(\theta_\chi\omega,\zeta_0,(n,\zeta_n) )$
is a partial sum of a stationary ergodic sequence.

%FC: By Condition L, there exist
%${\tilde\gamma}_1,{\tilde\gamma}_2>0$ such
% that $\muo_n(\partial\omega_n)\in
%({\tilde\gamma}_1,{\tilde\gamma}_2)$
By Condition \ref{ConditionL}, there exists ${\tilde\gamma}_1>0$ such
that $\mu^{\omega}_n(\partial\omega_n)\geq{\tilde\gamma}_1$
${\mathbb P}$-a.s.
Hence, since ${\mathbb P}$ is stationary and
%FC: $\kappa_{0,n}
$\kappa_{0,u} \in[0,1]$, we obtain for
$f \geq0$,
%e42 ###
%
\begin{eqnarray}\label{IP<IQ}
\langle{\mathrm E}^{\zeta}f(\theta_\chi\omega,\zeta_0) \rangle
_{{\mathbb P}} &=& \langle
\bar{\mathrm{E}}^{\zeta}f(\omega,\zeta_0) \rangle_{{\mathbb
P}}\nonumber\\
&=& \int_\Omega d{\mathbb P}\frac{1}{\mu^{\omega}_0(\partial\omega_0)}
\int_\Lambda d\mu^{\omega}_0(u)\,f(\omega,u)\nonumber\\[-8pt]\\[-8pt]
&\leq&\frac{1}{{\tilde\gamma}_1}\int_\Omega d{\mathbb P}
\int_\Lambda d\mu^{\omega}_0(u)\, \kappa_{0,u}^{-1}f(\omega,u)
\nonumber\\
&=& \frac{{\mathcal Z}}{{\tilde\gamma}_1}\langle f \rangle_{{\mathbb
Q}},\nonumber
\end{eqnarray}
which implies that
\[
%% \label{E_IP_f2<infty}
\mbox{if }f\in L^2({\mathbb Q}) \qquad\mbox{then }
\langle{\mathrm E}^{\zeta}f^2(\theta_\chi\omega,\zeta_0) \rangle
_{{\mathbb P}}<\infty.
\]

To proceed, we need to show that the random tube satisfies a
uniform local D\"{o}blin condition.
Denote ${\widetilde K}^{(n)} := K+K^{(2)}+\cdots+K^{(n)}$.
\begin{lmm}
\label{l_local_Doeblin}
Under Condition \ref{ConditionP}, there exist $N$ and ${\hat\gamma}>0$ such that
for all $x,y\in{\mathcal R}_{\omega}$ with $|(x-y)\cdot e|\leq2$ it
holds that
${\widetilde K}^{(N)}(x,y)\geq{\hat\gamma}$, ${\mathbb P}$-a.s.
\end{lmm}
\begin{pf*}{Proof of Lemma \protect\ref{l_local_Doeblin}}
Indeed, with the notation used in Condition \ref{ConditionP} and $n=N-1$,
\begin{eqnarray*}
{\widetilde K}^{(N)}(x,y) &\geq& K^{(n+1)}(x,y)\\
&\geq& \int_{B_1}K(x,z_1) \,d\nu^{\omega}(z_1) \\
&&{} \times\int_{B_2}K(z_1,z_2) \,d\nu^{\omega}(z_2)\cdots
\int_{B_n}K(z_{n-1},z_n)K(z_n,y) \,d\nu^{\omega}(z_n) \\
&\geq& \delta^n\varepsilon^{n+1}\\
&\geq& \delta^{N-1}\varepsilon^N.
\end{eqnarray*}
\upqed\end{pf*}

%%%%%%%%%%%%%%%%%%%%%%%%%%%%%
Next, we state some integrability and centering properties
which will be needed later.
\begin{lmm}
\label{l_resume1617}
We have
%e43 ###
%
\begin{eqnarray}
\label{EQ_corr2<infty}
\langle{\mathrm E}^{\zeta}\psi^2 (\omega,u,(\chi,\zeta_0) )
\rangle_{{\mathbb Q}} &<&
\infty,
\\
%
%e44 ###
%
\label{E_corr=0}
\langle{\mathrm E}^{\zeta}\psi(\theta_\chi\omega,\zeta
_0,(1,\zeta_1) ) \rangle
_{{\mathbb P}} &=&
\langle\bar{\mathrm{E}}^{\zeta}\psi(\omega,\zeta_0,(1,\zeta_1)
) \rangle_{{\mathbb P}} = 0.
\end{eqnarray}
\end{lmm}
\begin{pf*}{Proof of Lemma \protect\ref{l_resume1617}}
We start proving that
%e45 ###
%
\begin{equation}
\label{E_corr2<infty}
\langle{\mathrm E}^{\zeta}\psi^2 (\theta_\chi\omega,\zeta
_0,(1,\zeta_1) )
\rangle_{{\mathbb P}}
< \infty.
\end{equation}

We know that $\operatorname{grad}\psi\in L^2({\mathfrak N})$, that is,
$\langle(\operatorname{grad}\psi)^2 \rangle_{{K{\mathbb Q}}}<\infty$.
Analogously to the proof of the symmetry
relation (\ref{calculation_symmetry}), we obtain [note also that,
by (\ref{gradient=-diverg}), $\langle\operatorname{div}g \rangle
_{{\mathbb Q}}=0$
for all $g\in L^2({\mathfrak N})$]
\begin{eqnarray*}
\langle(\operatorname{grad}\psi)^2 \rangle_{{K{\mathbb Q}}} &=&
\biggl\langle\int_{{\mathsf V}
^+_{\omega,u}} \psi^2(\omega,u,y)K ((0,u),y ) \,d\nu^{\omega}(y)
\biggr\rangle
_{{\mathbb Q}} \\
&=& \biggl\langle\int_{{\mathsf V}^-_{\omega,u}} \psi^2(\omega,u,y)K
((0,u),y )\,
d\nu^{\omega}(y) \biggr\rangle_{{\mathbb Q}} \\
&=& \frac{1}{2} \bigl\langle{\mathtt E}_{\omega}\bigl(\psi^2(\omega,u,\xi
_1)\mid\xi
_0=(0,u)\bigr) \bigr\rangle_{{\mathbb Q}}.
\end{eqnarray*}
Then, using (\ref{cocycle}) we write
\[
\psi^2(\omega,u,\xi_n) \leq n\psi^2(\omega,u,\xi_1)
+ n\sum_{j=1}^{n-1}
\psi^2 \bigl(\theta_{\xi_j\cdot e}\omega,{\mathcal U}\xi_j,
\xi_{j+1}-(\xi_j\cdot e)e \bigr).
\]
Since ${\mathbb Q}$ is reversible for $G$, this implies that
for any $n$
%e46 ###
%
\begin{equation}
\label{corr_in_L2}
\biggl\langle\int_{\partial\omega}\psi^2(\omega,u,y) K^{(n)} ((0,u),y
) \,d\nu^{\omega}(y) \biggr\rangle_{{\mathbb Q}} < \infty,
\end{equation}
where $K^{(n)}(x,y)$ is the $n$-step transition density.

Let us define
\[
F_n^\omega= \{x\in\partial\omega\dvtx x\cdot e \in(n-1/2,n+1/2]
\}.
\]
Now we are going to use (\ref{corr_in_L2})
and Lemma \ref{l_local_Doeblin} to prove (\ref{E_corr2<infty}).
Note that, by Condition \ref{ConditionL}, there are positive constants
$C_1,C_2$ such that
\[
%% \label{bounded_F}
C_1 \leq\nu^{\omega}(F_1^\omega) \leq C_2,\qquad \mbox{${\mathbb P}$-a.s.}
\]
Using (\ref{cocycle}), we write on $\{\chi=0\}$
\begin{eqnarray*}
&&\psi^2 (\omega,\zeta_0,(1,\zeta_1) )\\
&&\qquad= \frac{1}{\nu^{\omega}(F_1^\omega)}\int_{F_1^\omega}
\psi^2 (\omega,\zeta_0,(1,\zeta_1) ) \,d\nu^{\omega}(y)\\
&&\qquad \leq\frac{2}{\nu^{\omega}(F_1^\omega)} \int_{F_1^\omega}
\bigl(\psi^2(\omega,\zeta_0,y)
+\psi^2(\theta_1\omega,\zeta_1,y-e) \bigr)\, d\nu^{\omega}(y)\\
&&\qquad\leq\frac{2}{{\hat\gamma}C_1} \biggl(\int_{\partial\omega}
{\widetilde K}^{(N)} ((0,\zeta_0),y )
\psi^2(\omega,\zeta_0,y) \,d\nu^{\omega}(y)\\
&&\qquad\quad\hspace*{26.6pt} {}+ \int_{\partial\omega}
{\widetilde K}^{(N)} ((1,\zeta_1),y )
\psi^2(\theta_1\omega,\zeta_1,y-e)\,
d\nu^{\omega}(y) \biggr).
\end{eqnarray*}
Using the stationarity of $\zeta$ under ${\mathbb P}$, we obtain that
\[
\langle{\mathrm E}^{\zeta}\psi^2 (\theta_\chi\omega,\zeta
_0,(1,\zeta_1) )
\rangle_{{\mathbb P}}
= \langle\bar{\mathrm{E}}^{\zeta}\psi^2 (\omega,\zeta_0,(1,\zeta
_1) ) \rangle_{{\mathbb P}},
\]
then, again by stationarity,
\begin{eqnarray*}
&&\biggl\langle\bar{\mathrm{E}}^{\zeta}\int_{\partial\omega}
{\widetilde K}^{(N)}
((1,\zeta_1),y ) \psi^2(\theta_1\omega,\zeta_1,y-e) \,d\nu^{\omega}(y)
\biggr\rangle_{{\mathbb P}} \\
&&\qquad= \biggl\langle\bar{\mathrm{E}}^{\zeta}\int_{\partial\omega}
{\widetilde K}^{(N)}
((0,\zeta_0),y ) \psi^2(\omega,\zeta_0,y) \,d\nu^{\omega}(y)
\biggr\rangle_{{\mathbb P}}.
\end{eqnarray*}
So, (\ref{E_corr2<infty}) follows from (\ref{IP<IQ})
and (\ref{corr_in_L2}).

Analogously, it is not difficult to prove that
(\ref{EQ_corr2<infty}) holds.
Indeed,
similarly to (\ref{IP<IQ}),
we have %% (recall the calculation \eqref{IP<IQ})
\begin{eqnarray*}
{\mathrm E}^{\zeta}\psi^2 (\omega,u,(\chi,\zeta_0) ) &=&
\int_{-1/2}^{1/2}d\alpha
\frac{1}{\mu^{\omega}_\alpha(\partial\omega_\alpha)}
\int_\Lambda d\mu^{\omega}_\alpha(v)\,\psi^2 (\omega,u,(\alpha,v)
)\\
&\leq&\frac{1}{{\tilde\gamma}_1}\int_{-1/2}^{1/2}d\alpha
\int_\Lambda d\mu^{\omega}_\alpha(v)\,\kappa_{\alpha,v}^{-1}
\psi^2 (\omega,u,(\alpha,v) )\\
&=& \frac{1}{{\tilde\gamma}_1}\int_{F_0^\omega}
\psi^2(\omega,u,y) \,d\nu^{\omega}(y),
\end{eqnarray*}
where we used (\ref{differentials}) in the last equality.
So, by Lemma \ref{l_local_Doeblin},
\[
{\mathrm E}^{\zeta}\psi^2 (\omega,u,(\chi,\zeta_0) ) \leq
\frac{1}{{\hat\gamma}{\tilde\gamma}_1}
\int_{\partial\omega}{\widetilde K}^{(N)} ((0,u),y )
\psi^2(\omega,u,y) \,d\nu^{\omega}(y)
\]
and (\ref{EQ_corr2<infty}) follows from (\ref{corr_in_L2}).

Finally, let us prove (\ref{E_corr=0}).
The first equality follows from the stationarity of ${\mathbb P}$.
Then, since $\operatorname{grad}\psi\in L^2_\nabla({\mathfrak N})$,
there is a sequence
of functions $f_n\in L^2({\mathfrak S})$ such that $\nabla f_n\to
\operatorname{grad}\psi$
in the sense of the $L^2({\mathfrak N})$-convergence.
Note that,
in fact, when proving (\ref{E_corr2<infty}), we proved that for
any function $g\in{\mathcal H}({\mathfrak M})$ such that
$\operatorname{grad}g \in L^2_\nabla({\mathfrak S})$,
we have for some $C_3>0$,
\[
\langle{\mathrm E}^{\zeta}g^2 (\theta_\chi\omega,\zeta_0,(1,\zeta
_1) ) \rangle
_{{\mathbb P}} <
C_3\langle(\operatorname{grad}g)^2 \rangle_{{K{\mathbb Q}}}.
\]
Then, (\ref{E_corr=0}) follows from the
above fact applied to $g$
assembled from shifts of $\operatorname{grad}\psi- \nabla f_n$,
since then we can then write
\[
\langle\psi(\theta_\chi\omega,\zeta_0,(1,\zeta_1) ) \rangle
_{{\mathbb P}}
= \lim_{n\to\infty} [\langle f_n(\theta_{\chi+1} \omega,\zeta_1)
\rangle_{{\mathbb P}} - \langle f_n(\theta_{\chi} \omega,\zeta_0)
\rangle
_{{\mathbb P}} ]=0
\]
by the stationarity of ${\mathbb P}$.
\end{pf*}

%s4 ###
\section{Proofs of the main results}
\label{s_proofs}

%s4.1 ###
\subsection{\texorpdfstring{Proof of Theorem
\protect\ref{t_q_invar_princ}}{Proof of Theorem 2.1}}
\label{s_proof_inv_pr}
In this section, we apply the machinery of
Section~\ref{s_environment} in order to prove the invariance
principle
for the (discrete time) motion of a single particle.
\begin{pf*}{Proof of Theorem \protect\ref{t_q_invar_princ}}
Denote
\[
\Theta_n = \xi_n \cdot e +
\psi(\theta_\chi\omega,\zeta_0,\xi_n-\chi e).
\]
Observe that by (\ref{psi_is_corrector}), $\Theta$ is a martingale
under the quenched law ${\mathtt P}_{\omega}$.
By shift-covariance (\ref{shift-covariance}) the increments of
$\Theta_n$ do not depend of $\chi$ and $\zeta$. With the notation
\[
h(\omega,u)= {\mathtt E}_{\omega}[ (\Theta_1-\Theta_0)^2 \mid\xi_0=(0,u)],
\]
the bracket of the martingale $\Theta_n$ is given by
\[
\langle\Theta\rangle_n = \sum_{i=0}^{n-1}
h( \theta_{\xi_i\cdot e}\omega, {\mathcal U}\xi_{i}).
\]
By the ergodic theorem,
%e47 ###
%
\begin{equation}
\label{def:sigma}
\frac{1}{n} \langle\Theta\rangle_n \longrightarrow
\sigma^2 \stackrel{\mathrm{def}}{=} \langle h(\omega,u) \rangle
_{{\mathbb Q}},
\end{equation}
a.s. as $n \to\infty$.
Clearly, $\sigma^2 \in(0,\infty)$. Moreover, for all $\varepsilon>0$,
%e48 ###
%
\begin{equation}
\label{eq:cnrs1}
\sum_{i=0}^{n-1} {\mathtt P}_{\omega}[ |\Theta_{i+1}-\Theta_{i}|
\geq
\varepsilon n^{1/2} \mid\xi_i ] \to0
\end{equation}
for ${\mathbb P}$-a.e. $\omega$ and ${\mathtt P}_{\omega}$-a.e. path.
To show this,
define for any $a>0$ and all $n\geq1$,
\[
h^{(a)}_n(\omega) = {\mathtt E}_{\omega}\bigl((\Theta_n-\Theta_{n-1})^2
{\mathbb I}{\{|\Theta_n-\Theta_{n-1}|\geq a \}} \mid\xi_{n-1} \bigr).
\]
Using the ergodicity of the process of the environment
viewed from the particle, we obtain
\[
\frac{1}{n} \sum_{i=1}^n h^{(a)}_i
\longrightarrow
\bigl\langle{\mathtt E}_{\omega}\bigl((\Theta_1-\Theta_0)^2{\mathbb I}{\{
|\Theta _1-\Theta_0|\geq a\}}
\mid\xi_0 = (0,u) \bigr) \bigr\rangle_{{\mathbb Q}}
\]
as $n \to\infty$
for ${\mathbb P}$-almost all $\omega$ and ${\mathtt P}_{\omega}$-almost
all trajectories of the walk. Note that, when $a$ is replaced
by $\varepsilon n^{1/2}$,
the left-hand side is, by
Bienaym\'{e}--Chebyshev inequality, an upper bound of the left-hand
side of (\ref{eq:cnrs1}) multiplied by $\varepsilon^2$.
Hence (\ref{eq:cnrs1}) follows by taking $a$ arbitrarily large.

% Define for any $\eps>0$ and all $n\geq1$
% \[
% f^{(\eps)}_n(\omega) = \Eo((\Theta_n-\Theta_{n-1})^2
% \vienas{|\Theta_n-\Theta_{n-1}|\geq\eps} ).
% \]
% Using the ergodicity of the process of the environment
% viewed from the particle, we obtain
% \[
% \frac{1}{n} \sum_{i=1}^n f^{(\eps)}_i
% \to\exq{\Eo((\xi_1\cdot e)^2\vienas{|\xi_1\cdot e|\geq\eps}
% \mid\xi_0 = (0,u) )}
% \]
% for $\IP$-almost all $\omega$ and $\Po$-almost
% all trajectories of the Knudsen random walk.
% \dots So,
Combining (\ref{def:sigma}) and (\ref{eq:cnrs1}),
we can apply the central limit theorem for martingales
(cf., e.g., Theorem 7.7.4 of \cite{D})
to show that
%e49 ###
%
\begin{equation}
\label{conv_to_BM}
%% n^{-1/2}\Theta_{[n\cdot]} \to\sigma B(\cdot)
%% { in law as $n\to\infty$},
n^{-1/2}\Theta_{[n\cdot]}
\stackrel{\mathrm{law}}{\longrightarrow}
\sigma B(\cdot)\qquad
\mbox{as $n\to\infty$},
\end{equation}
where $B(\cdot)$ is the Brownian motion.

Then the idea is the following:
using (\ref{E_corr=0}) and the ergodic theorem, we
obtain that the
%FC: corrector $\psi$ behaves sublinearly,
corrector $\psi(\omega,u,x)$ behaves sublinearly in $x$
which implies the convergence of $n^{-1/2}\xi_{[nt]}\cdot e$.
More precisely, we can write, with $m_j:=[1/2+\xi_j\cdot e]$ and
using (\ref{shift-covariance}),
%e50 ###
%
\begin{eqnarray}\label{decomp_mart}
\frac{\Theta_{[nt]}}{n^{1/2}}
&=& \frac{\xi_{[nt]}\cdot e +\psi(\theta_\chi\omega,\zeta_0,
(m_{[nt]},\zeta_{m_{[nt]}}) )}
{n^{1/2}}\nonumber\\[-8pt]\\[-8pt]
&&{} - \frac{\psi(
\theta_{\xi_{[nt]}\cdot e}\omega,{\mathcal U}\xi_{[nt]},
%% (m_{[nt]},\zeta_{m_{[nt]}}) )}{n^{1/2}}.
(\chi+m_{[nt]}-\xi_{[nt]}\cdot e,\zeta_{m_{[nt]}})
)}{n^{1/2}}.\nonumber
\end{eqnarray}
Let us prove that the second term in the right-hand side
converges to $0$ in ${\mathtt P}_{\omega}$-probability for ${\mathbb P}$-almost
all $\omega$ and almost all $(\chi,\zeta)$.
Suppose, for the sake of simplicity, that $t=1$.
Then, by the stationarity of the process $ ((\chi+n,\zeta_n),
n\in{\mathbb Z})$ and (\ref{IQ_is_invariant}) together
with (\ref{EQ_corr2<infty}), we have for all $i\geq0$,
\begin{eqnarray*}
&&\bigl\langle{\mathrm E}^{\zeta}{\mathtt E}_{\omega}\bigl[\psi^2
\bigl(\theta_{\xi_i\cdot e}\omega, {\mathcal U}
\xi_i, (\chi+m_i-\xi_i\cdot e,\zeta_{m_i}) \bigr)\mid\xi_0=(0,u) \bigr]
\bigr\rangle_{{\mathbb Q}}
\\
&&\qquad=\langle{\mathrm E}^{\zeta}\psi^2 (\omega,u,(\chi,\zeta_0) )
\rangle_{{\mathbb Q}}\\
&&\qquad<\infty,
\end{eqnarray*}
%
%%% This is wrong:
% \begin{eqnarray*}
% \lefteqn{ \exip{\Ez
% \Eo\psi^2 (\theta_{\chi+\xi_i\cdote}\omega,
% \UU\xi_i,(m_i,\zeta_{m_i}) )}} \\
% & = \exip{\Ez\Eo\psi^2 (
% \theta_{\chi+\xi_{i+1}\cdot e}\omega,
% \UU\xi_{i+1},(m_{i+1},\zeta_{m_{i+1}}) )},
% \end{eqnarray*}
so, by the ergodic theorem,
\begin{eqnarray*}
&&\frac{1}{n}\sum_{i=1}^n
{\mathtt E}_{\omega}\bigl(\psi^2 (\theta_{\chi+\xi_i\cdot e}\omega
,{\mathcal U}\xi_i,
(m_i,\zeta_{m_i}) )\\
&&\hspace*{17.6pt}\qquad{} - \psi^2 (\theta_{\chi+\xi_{i-1}\cdot e}\omega,
{\mathcal U}\xi_{i-1}, (m_{i-1},\zeta_{m_{i-1}}) ) \bigr) \to0
\end{eqnarray*}
as $n\to\infty$ which implies that
%e51 ###
%
\begin{equation}
\label{2nd_to_0}
\frac{1}{n}{\mathtt E}_{\omega}\psi^2 (\theta_{\chi+\xi_n\cdot
e}\omega,
{\mathcal U}\xi_n,(m_n,\zeta_{m_n}) )
\to0
\end{equation}
for ${\mathbb P}$-almost all $\omega$ and almost all $(\chi,\zeta)$.
Now, let us prove that the limit of the first term
in the right-hand side of (\ref{decomp_mart})
is the same as the limit of $n^{-1/2}\xi_{[nt]}\cdot e$; for this,
we have to prove that
%e52 ###
%
\begin{equation}
\label{psi_to_0}
\frac{\psi(\theta_\chi\omega,\zeta_0,(m_{[nt]},
\zeta_{m_{[nt]}}) )}{n^{1/2}}
\to0 \qquad\mbox{as }n\to\infty, \mbox{ in ${\mathtt P}_{\omega}$-probability}.
\end{equation}
Using (\ref{decomp_corr}), (\ref{E_corr=0}), and the ergodic
theorem, we obtain that for ${\mathbb P}$-almost all $\omega$
$m^{-1}\psi(\theta_\chi\omega,\zeta_0,(m,\zeta_m) )\to0$
for almost all $(\chi,\zeta)$, as $|m| \to\infty$.
This means that, for any $\varepsilon>0$ there exists $H$
(depending on $\omega,\zeta,\chi$)
such that
%e53 ###
%
\begin{equation}
\label{psi_sublinear}
|\psi(\theta_\chi\omega,\zeta_0,(m,\zeta_m) ) |
\leq H+\varepsilon|m|.
\end{equation}
Denote
\[
\Psi_j = \xi_j\cdot e +
\psi(\theta_\chi\omega,\zeta_0,(m_j,\zeta_{m_j}) ).
\]
From (\ref{psi_sublinear}) we see that
\begin{eqnarray*}
|\psi(\theta_\chi\omega,\zeta_0,(m_j,\zeta_{m_j}) )
| &\leq& H+\varepsilon|m_j|\\
&\leq& H+\frac{\varepsilon}{2}+\varepsilon|\xi_j\cdot e|\\
&\leq& H+\frac{\varepsilon}{2}+\varepsilon\bigl(|\Psi_j|
+ |\psi(\theta_\chi\omega,
\zeta_0,(m_j,\zeta_{m_j}) ) | \bigr),
\end{eqnarray*}
so for $\varepsilon<1/2$ we obtain
\[
|\psi(\theta_\chi\omega,\zeta_0,(m_j,\zeta_{m_j}) )
| \leq2H+\varepsilon+2\varepsilon|\Psi_j|.
\]
Using (\ref{conv_to_BM}) and (\ref{2nd_to_0})
in (\ref{decomp_mart}), we obtain
\[
\max_{j\leq n} \frac{|\Psi_j|}{n^{1/2}}
\stackrel{\mathrm{law}}{\longrightarrow} \sigma
{\max_{s\in[0,1]}}|B(s)|.
\]
So by the portmanteau theorem (cf. Theorem 2.1(iii)
of \cite{Bil}),
\[
\limsup_{n\to\infty} {\mathtt P}_{\omega}\Bigl[{\max_{j\leq n}}
|\psi(\theta_\chi\omega,\zeta_0,(m_j,\zeta_{m_j}) )
|\geq an^{1/2} \Bigr]
\leq P \biggl[{\max_{s\in[0,1]}}|B(s)|\geq\frac{a \sigma}{2\varepsilon}
\biggr],
\]
which converges to $0$ for any $a$ as $\varepsilon\to0$.
This concludes the proof of Theorem~\ref{t_q_invar_princ}.
%\rightqed
\end{pf*}

%s4.2 ###
\subsection{On the finiteness of the second moment}
\label{s_b_finite}
In this section, we prove the results which concern the
finiteness of $\langle b \rangle_{{\mathbb Q}}$.
First, we present a (quite elementary)
proof of Proposition \ref{pr_suff_2nd_moment} in the case $d\geq4$.
\begin{pf*}{Proof of Proposition \protect\ref{pr_suff_2nd_moment}
\textup{(case $d\geq4$)}}
First of all, note that
\[
|\{s\in{\mathbb S}^{d-1}\dvtx x+hs\in{\mathbb R}\times\Lambda\}| = O\bigl(h^{-(d-1)}\bigr)
\qquad\mbox{as } h\to\infty,
\]
uniformly in $x\in{\mathbb R}\times\Lambda$.
So, since $\omega\subset{\mathbb R}\times\Lambda$,
there is a constant $C_1>0$, depending only
on ${\widehat M}=\operatorname{diam}(\Lambda)/2$ and the dimension,
such that for ${\mathbb P}$-almost all~$\omega$
%e54 ###
%
\begin{equation}
\label{bound_d>=4}
{\mathtt P}_{\omega}[|(\xi_1-\xi_0)\cdot e| > h \mid\xi_0=x] \leq
C_1 h^{-(d-1)}
\end{equation}
for all $x\in\partial\omega$, $h\geq1$.
Inequality (\ref{bound_d>=4}) immediately implies that $b$
is uniformly bounded for $d\geq4$.
\end{pf*}

Unfortunately, the above proof does not work in the case $d=3$.
To treat this case, we need some results concerning induced chords
which in some sense generalize
Theorems 2.7 and 2.8 of \cite{CPSV1}. So the rest of this section
is organized as follows. After introducing some notation, we prove
Proposition \ref{p_cosine_int} which is a generalization of the
result about the induced chord in a convex subdomain (Theorem 2.7
of \cite{CPSV1}). This will allow us to prove
Proposition \ref{pr_infinite_d=2}.
Then, using Theorem 2.8 of \cite{CPSV1} (the result about induced
chords in a general subdomain) we prove
Proposition~\ref{p_induced_chords}---a property of random
chords induced in a smaller random tube by
a random chord in a bigger random tube.
This last result will allow us to
prove Proposition \ref{pr_suff_2nd_moment}.

\begin{figure}

\includegraphics{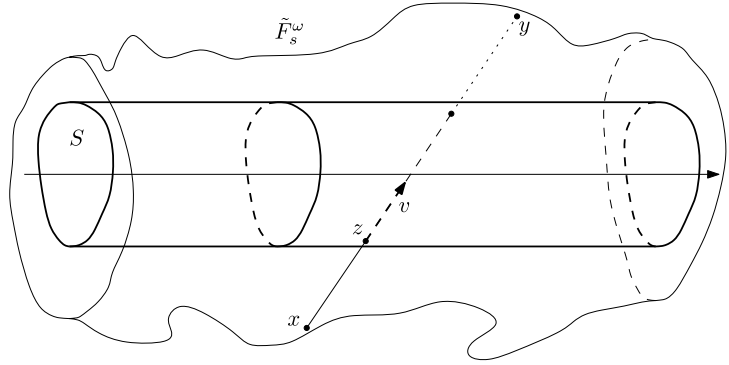}

\caption{On the definition of ${\mathsf L}$ and ${\mathsf Y}$: we have
${\mathsf L}(x,y)={\mathcal U}z$ and ${\mathsf
Y}(x,y)=U^{-1}_{{\mathcal U}z}v$.}
\label{f_cylinder}
\end{figure}

Let $S\subset\Lambda$ be an open convex set, and denote
by ${\widehat S}= {\mathbb R}\times S$ the straight cylinder generated
by $S$.
Assuming that ${\widehat S}\subset\omega$, we let ${\mathcal I}$ be
the event that the
trajectory of the first jump (i.e., from $\xi_0$ to $\xi_1$)
intersects ${\widehat S}$:
\[
{\mathcal I}= \{\mbox{there exists }t\in[0,1]\mbox{ such that }\xi_0
+(\xi_1-\xi_0)t\in{\widehat S}\}.
\]
For any $u\in\partial S$ such that $\partial S$ is differentiable
in $u$, define ${\hat{\mathbf n}}_u$ to be the normal vector
with respect to $\partial{\widehat S}$ at the point $(0,u)$; clearly,
we have ${\hat{\mathbf n}}_u\cdot e=0$ (if $\partial S$ is not differentiable
in $u$, define ${\hat{\mathbf n}}_u$ arbitrarily). Fix some family
$(U_v, v\in
\partial S)$ of unitary linear operators with the property $U_v e =
{\hat{\mathbf n}}_v$ for all $v\in\partial S$. Now, on the event
${\mathcal I}$ we may
define the conditional law of intersection of $\partial{\widehat S}$. Namely,
for $x,y\in\partial\omega$, let
%e55 ###
%
\begin{equation}
\label{df_t_xy}
t_{x,y} = \inf\{t\in[0,1]\dvtx x+(y-x)t \in\partial{\widehat S}\}
\end{equation}
with the convention $\inf\varnothing=\infty$.
Then, we define the (projected) location of the crossing
of $\partial{\widehat S}$
by
\[
{\mathsf L}(x,y)= \cases{
{\mathcal U}\bigl(x+(y-x)t_{x,y}\bigr), &\quad if $t_{x,y} \in[0,1]$,\cr
\infty, &\quad otherwise,}
\]
and the relative direction of the crossing by
\[
{\mathsf Y}(x,y)= \cases{
\displaystyle U_{{\mathsf L}(x,y)}^{-1}\frac{y-x}{\|y-x\|}, &\quad if $t_{x,y}
\in[0,1]$,\vspace*{2pt}\cr
0, &\quad otherwise,}
\]
(see Figure \ref{f_cylinder}).

%f3 ###
%

Here, in the case when there is no intersection, for formal
reasons we still assign values for ${\mathsf L}$ and ${\mathsf Y}$;
note, however,
that in the case $t_{x,y} \in[0,1]$, we have
${\mathsf L}(x,y)\in\partial S$ and ${\mathsf Y}(x,y)\in{\mathbb S}_e$.

Before proving Proposition \ref{pr_infinite_d=2}, we obtain a
remarkable fact which is
closely related to the invariance properties of random chords (cf.
Theorems 2.7 and 2.8 of \cite{CPSV1}).
We have that, conditioned on ${\mathcal I}$, the annealed law of
the pair of random variables
$({\mathsf L}(\xi_0,\xi_1),{\mathsf Y}(\xi_0,\xi_1))$ can be described
as follows: ${\mathsf L}(\xi_0,\xi_1)$ and ${\mathsf Y}(\xi_0,\xi
_1)$ are
independent, ${\mathsf L}(\xi_0,\xi_1)$ is uniform
on $\partial S$ and ${\mathsf Y}(\xi_0,\xi_1)$ has the cosine distribution.
More precisely, we formulate and prove the following result.
\begin{prop}
\label{p_cosine_int}
Let $d \geq2$.
It holds that $\langle{\mathtt P}_{\omega}[{\mathcal I}] \rangle
_{{\mathbb Q}}={|\partial S|}/{{\mathcal Z}
}$. Moreover,
for any measurable $B_1\subset\partial S, B_2\subset{\mathbb S}_e$
we have
%e56 ###
%
\begin{eqnarray}
\label{eq_cosine_int}
% \frac{1}{\exq{\Po[\I]}}
&&\bigl\langle{\mathtt E}_{\omega}\bigl({\mathbb I}{\{{\mathsf L}(\xi_0,\xi
_1)\in B_1, {\mathsf Y}(\xi_0,\xi_1)\in B_2\}}\mid\xi_0=(0,u) \bigr)
\bigr\rangle_{{\mathbb Q}}\nonumber\\[-8pt]\\[-8pt]
&&\qquad= \frac{|\partial S|}{{\mathcal Z}} \frac{|B_1|}{|\partial S|}
\gamma_d \int_{B_2} h\cdot e \,dh.\nonumber
\end{eqnarray}
\end{prop}
\begin{pf}
First, we prove (\ref{eq_cosine_int}).
Define ${\tilde F}^\omega_s = \{x\in\partial\omega:
x\cdot e\in[-s,s]\}$ for
$s>0$. By the translation invariance and (\ref{differentials}),
we have
%e57 ###
%
\begin{eqnarray}\label{conta_inv*}\quad
&&\bigl\langle{\mathtt E}_{\omega}\bigl({\mathbb I}{\{{\mathsf L}(\xi
_0,\xi _1)\in B_1, {\mathsf Y}(\xi_0,\xi _1)\in B_2\}}\mid\xi
_0=(0,u) \bigr) \bigr\rangle_{{\mathbb Q}}
\nonumber\\
&&\qquad= \frac{1}{{\mathcal Z}}\int_\Omega d{\mathbb P}
\int_\Lambda d\mu^{\omega}_0(u)\,\kappa_{0,u}^{-1}\nonumber\\
&&\qquad\quad{} \times\int_{\partial
\omega}
d\nu^{\omega}(y)\,K ((0,u),y ){\mathbb I}{\{{\mathsf L}((0,u),y) \in B_1, {\mathsf
Y}((0,u),y)\in B_2\}}\nonumber\\
&&\qquad= \frac{1}{2s{\mathcal Z}}\int_\Omega d{\mathbb P}\int_{-s}^s ds
\int_\Lambda d\mu^{\omega}_s(u)\,\kappa_{s,u}^{-1}\\
&&\qquad\quad{} \times\int_{\partial
\omega}
d\nu^{\omega}(y)\,K ((s,u),y )
{\mathbb I}{\{{\mathsf L}((s,u),y) \in B_1, {\mathsf
Y}((s,u),y)\in B_2\}}\nonumber\\
&&\qquad= \frac{1}{2s{\mathcal Z}}\int_\Omega d{\mathbb P}\int_{{\tilde
F}^\omega_s}
d\nu^{\omega}(x)\nonumber\\
&&\qquad\quad{}\times\int_{\partial\omega}d\nu^{\omega}(y)\, {\mathbb
I}{\{{\mathsf L}(x,y)\in B_1, {\mathsf Y}(x,y)\in
B_2\}}K(x,y).\nonumber
\end{eqnarray}
Define the domain ${\mathcal D}_s^\omega$ by
\[
{\mathcal D}_s^\omega= \{x\in\omega\dvtx x\cdot e \in[-s,s]\}
\]
and note that $\partial{\mathcal D}_s^\omega= {\tilde F}^\omega
_s\cup
(\{-s\}\times\omega_{-s})\cup(\{s\}\times\omega_{s})$.
For $x,y\in\partial{\mathcal D}_s^\omega$ let ${\hat K}(x,y)$ be defined
as in (\ref{def_K}), but with ${\mathcal D}_s^\omega$ instead
of $\omega$. Note that ${\hat K}(x,y) = K(x,y)$ when
$x,y\in{\tilde F}^\omega_s$.

Next, we show that the random chord in $\omega$ with the first
point in ${\tilde F}^\omega_s$ has roughly the same
law as the random chord in ${\mathcal D}_s^\omega$: for any
$\varepsilon>0$
there exists $s_0$ such that for all $s\geq s_0$ [with some abuse of
notation, we still denote by $\nu^{\omega}(\partial{\mathcal
D}_s^\omega)$ the
$(d-1)$-dimensional Hausdorff measure of $\partial{\mathcal
D}_s^\omega$]
%e58 ###
%
\begin{eqnarray}\label{conta_inv**}\hspace*{20pt}
&& \biggl| \frac{1}{\nu^{\omega}({\tilde F}^\omega_s)}
\int_{{\tilde F}^\omega_s} d\nu^{\omega}(x)\int_{\partial\omega
}d\nu^{\omega}(y)\,
{\mathbb I}{\{{\mathsf L}(x,y)\in B_1, {\mathsf Y}(x,y)\in B_2\}}K(x,y)
\nonumber\\
&&\qquad{} - \frac{1}{\nu^{\omega}(\partial{\mathcal D}_s^\omega)}
\int_{(\partial{\mathcal D}_s^\omega)^2}
d\nu^{\omega}(x)\,d\nu^{\omega}(y)\\
&&\,\hspace*{130.8pt}{}\times{\mathbb I}{\{{\mathsf L}(x,y)\in
B_1, {\mathsf Y}(x,y)\in B_2\}}
{\hat K}(x,y) \biggr| < \varepsilon\nonumber
\end{eqnarray}
[in the second term, we suppose that ${\mathsf L}(x,y)=\infty,
{\mathsf Y}(x,y)=0$ when $x\in(\{-s\}\times S)\cup(\{s\}\times S)$].
Indeed, we have
%e59 ###
%
\begin{equation}
\label{poverhnosti}
\nu^{\omega}({\tilde F}^\omega_s) \leq\nu^{\omega}(\partial
{\mathcal D}_s^\omega)
\leq\nu^{\omega}({\tilde F}^\omega_s)+2|\Lambda|
\end{equation}
and, by Condition \ref{ConditionL}, there exists $C_4>0$ such that
%e60 ###
%
\begin{equation}
\label{poverhnost'_Fs}
\nu^{\omega}({\tilde F}^\omega_s) \geq C_4 s, \qquad{\mathbb P}\mbox{-a.s.}
\end{equation}
Also, since $\omega\subset{\mathbb R}\times\Lambda$, for any
$\varepsilon>0$
there exists $C_5>0$ such that
for all $x\in\partial\omega$
%e61 ###
%
\begin{equation}
\label{chord_not_too_long}
\int_{\{y\in\partial\omega\dvtx |(x-y)\cdot e|>C_5\}}
K(x,y) \,d\nu^{\omega}(y) < \varepsilon,\qquad
{\mathbb P}\mbox{-a.s.}
\end{equation}
Now, (\ref{conta_inv**}) follows from
(\ref{poverhnosti})--(\ref{chord_not_too_long}) and
a coupling argument: choose the first point uniformly
on $\partial{\mathcal D}_s^\omega$; with
big probability, it will fall on ${\tilde F}^\omega_{s-C_5}$
(and so it can be used
as the first point of the random chord in $\partial\omega$).
Then, choose the second point according
to the cosine law; by (\ref{chord_not_too_long}), with big
probability it will belong to ${\tilde F}^\omega_s$, and so the
probability that the coupling is successful converges to $1$
as $s\to\infty$.

Then, recall Theorem 2.7 from \cite{CPSV1}: in a finite domain,
the induced random chord in a convex subdomain has the same
uniform${}\times{}$cosine law. So
\begin{eqnarray*}
&&\frac{1}{\nu^{\omega}(\partial{\mathcal D}_s^\omega)}
\int_{\partial{\mathcal D}_s^\omega\times\partial{\mathcal
D}_s^\omega}
d\nu^{\omega}(x)\,d\nu^{\omega}(y)\,{\mathbb I}{\{{\mathsf L}(x,y)\in
B_1, {\mathsf Y}(x,y)\in B_2\}}
{\hat K}(x,y)
\\
&&\qquad = {\mathtt P}_{\omega}[{\mathcal I}_s] \frac{|B_1|}{|{\mathbb
S}_e|} \gamma_d \int_{B_2}
h\cdot e \,dh,
\end{eqnarray*}
where ${\mathcal I}_s$ is the event that the random chord of
$\partial{\mathcal D}_s^\omega$
crosses the set $[-s,s]\times\partial S$.
By formula (47) of \cite{CPSV1}
[see also formula (17) in Theorem 2.8 there], we have
%e62 ###
%
\begin{equation}
\label{ver_peresech}
{\mathtt P}_{\omega}[{\mathcal I}_s] = \frac{2s|\partial
S|}{|\partial{\mathcal D}_s^\omega|}
= \frac{2s|\partial S|}{\nu^{\omega}({\tilde F}^\omega_s)
+|\omega_{-s}|+|\omega_s|}.
\end{equation}
Since, by the ergodic theorem, $|{\tilde F}^\omega_s|/(2s)\to
{\mathcal Z}$
as $s\to\infty$, (\ref{ver_peresech}) implies that
${\mathtt P}_{\omega}[{\mathcal I}_s]\to|\partial S|/{\mathcal Z}$ as
$s\to\infty$.
We obtain (\ref{eq_cosine_int}) using (\ref{conta_inv*})
and (\ref{conta_inv**}), and sending $s$ to $\infty$.

Finally, the fact that $\langle{\mathtt P}_{\omega}[{\mathcal I}]
\rangle_{{\mathbb Q}}={|\partial
S|}/{{\mathcal Z}}$
follows from (\ref{eq_cosine_int})
(take $B_1=\partial S, B_2={\mathbb S}_e$).
\end{pf}

Now, using Proposition \ref{p_cosine_int}, it is
straightforward to obtain Proposition \ref{pr_infinite_d=2}.
\begin{pf*}{Proof of Proposition \protect\ref{pr_infinite_d=2}}
Suppose that $\omega$ contains an infinite straight
cylinder ${\widehat S}$ (more precisely, a strip, since we
are considering the case $d=2$) of height $r>0$, ${\mathbb P}$-a.s.
Keep the notation $t_{x,y}$ from (\ref{df_t_xy}), and define also
\[
t'_{x,y} = \sup\{t\in[0,1]\dvtx x+(y-x)t \in\partial{\widehat S}\}.
\]
On the event ${\mathcal I}$, define the random points $\Upsilon_0,
\Upsilon_1\in\partial{\widehat S}$
by
\begin{eqnarray*}
\Upsilon_0 &=& \xi_0+(\xi_1-\xi_0)t_{x,y},\\
\Upsilon_1 &=& \xi_0+(\xi_1-\xi_0)t'_{x,y},
\end{eqnarray*}
so that $(\Upsilon_0, \Upsilon_1)$ is the random chord of ${\widehat S}$
induced by the first crossing. On ${\mathcal I}^c$, define
$\Upsilon_0=\Upsilon_1=0$. By Proposition \ref{p_cosine_int},
conditioned on ${\mathcal I}$, the random chord $(\Upsilon_0, \Upsilon
_1)$ has
the cosine law, that is, the density of a direction is proportional to
the cosine of the angle between this direction and the normal vector
(which, in this case, is perpendicular to $e$).
Let $P[\cdot]:=\frac{{\mathcal Z}}{2}\langle{\mathtt P}_{\omega
}[\cdot{\mathbb I}{\{\mathcal I\}}] \rangle_{{\mathbb Q}
}$ be the annealed
probability conditioned on the intersection;
since $d=2$ and $S$ is a bounded interval, $|\partial S|=2$.
With $\eta:=(\xi_1-\xi_0)/\|\xi_1-\xi_0\|$ and ${\hat{\mathbf
n}}$ the inner
normal vector to the cylinder at $\Upsilon_0$,
%f4 ###
%
\begin{figure}

\includegraphics{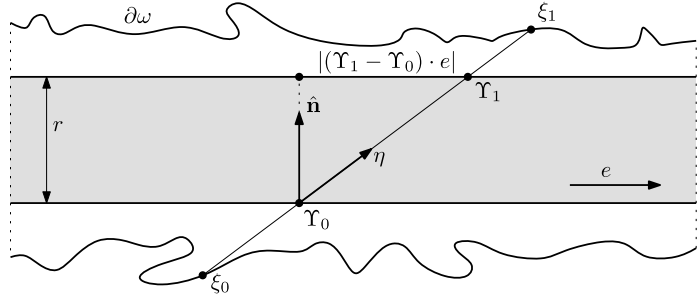}

\caption{($d=2$) Computing the distribution of
$|(\Upsilon_1-\Upsilon_0)\cdot e|$.}
\label{f_2dim}
\end{figure}
we have (see Figure \ref{f_2dim})
\begin{eqnarray*}
P [|(\Upsilon_1-\Upsilon_0)\cdot e| > x ] &=&
P \biggl[\eta\cdot{\hat{\mathbf n}}< \frac{r}{\sqrt{r^2+x^2}} \biggr]\\
&=& \int_{\arccos\frac{r}{\sqrt{r^2+x^2}}}^{\pi/2}\cos z \,dz\\
&=& 1-\frac{x}{\sqrt{r^2+x^2}},
\end{eqnarray*}
so the conditional density of the random variable
$|(\Upsilon_1-\Upsilon_0)\cdot e|$ is $f(x) =
\frac{r^2}{(r^2+x^2)^{3/2}}$ on
${\mathbb R}^+$.
Then we have
%FC: added explanations
% \exq{b}\geq\exq{\Eo|(\Upsilon_1-\Upsilon_0)\cdot e|^2}
% \times\exq{\Po[\I]}
% = \frac{2}{\ZZ} \intl_0^{+\infty} x^2
% \frac{r^2}{(r^2+x^2)^{3/2}} dx = +\infty,
%
\begin{eqnarray*}
\langle b \rangle_{{\mathbb Q}} &=&
\bigl\langle{\mathtt E}_{\omega}\bigl(|(\xi_1-\xi_0)\cdot e|^2 [{\mathbb I}{\{
{\mathcal I}\}}+{\mathbb I}{\{{\mathcal I}^c\}}] \mid\xi
_0=(0,u) \bigr) \bigr\rangle_{{\mathbb Q}}\\
&\geq&
\bigl\langle{\mathtt E}_{\omega}\bigl(|(\xi_1-\xi_0)\cdot e|^2 {\mathbb I}{\{
{\mathcal I}\}} \mid\xi_0=(0,u) \bigr)
\bigr\rangle_{{\mathbb Q}}\\
&\geq& \langle{\mathtt E}_{\omega}|(\Upsilon_1-\Upsilon_0)\cdot
e|^2 \rangle_{{\mathbb Q}
} \times
\langle{\mathtt P}_{\omega}[{\mathcal I}] \rangle_{{\mathbb Q}}\\
&=& \frac{2}{{\mathcal Z}} \int_0^{+\infty} x^2
\frac{r^2}{(r^2+x^2)^{3/2}} \,dx = +\infty,
\end{eqnarray*}
which concludes the proof of Proposition \ref{pr_infinite_d=2}.
\end{pf*}

\begin{figure}

\includegraphics{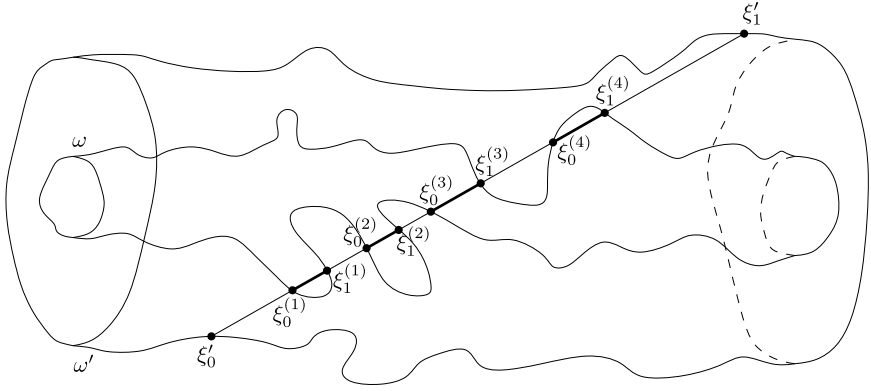}

\caption{Random chords induced in a random tube $\omega$ by a
random chord in a
random tube $\omega'$ (in this particular case, we have $\iota=4$).}
\label{f_many_chords}
\end{figure}

Let us observe that if a stationary ergodic random tube is
almost surely convex,
then necessarily it has the form ${\mathbb R}\times S$ for some convex
(and nonrandom) set $S\subset\Lambda$. This shows that
Proposition \ref{p_cosine_int} is indeed a generalization of
Theorem~2.7 of \cite{CPSV1}. Now our goal is to obtain an analogue
of a more general Theorem 2.8 of~\cite{CPSV1}. For this we consider
a pair of stationary ergodic random tubes
$(\omega,\omega')\in\Omega^2$,
let $\widetilde{{\mathbb P}}$ be their joint law and ${\mathbb
P},{\mathbb P}'$ be
the corresponding marginals. Suppose also that $\omega$ is
contained in $\omega'$ $\widetilde{{\mathbb P}}$-a.s.
We keep the notation such as $\kappa_x, K(x,y), \ldots$
for $x,y \in\partial\omega'$ as well, when it creates no
confusion; for the measures $\mu$ and $\nu$ we usually indicate in
the upper index whether they refer to $\omega$ or $\omega'$.
Denote also ${\mathcal Z}' = \int_{\Omega} d{\mathbb P}'\int_\Lambda
\kappa^{-1}_{0,u} \,d\mu_0^{\omega'}(u)$.
If $(\xi'_0,\xi'_1)$ is a chord in $\omega'$, we denote by
$(\xi_0^{(1)},\xi_1^{(1)}),
\ldots, (\xi_0^{(\iota)},\xi_1^{(\iota)})$ the induced random
chords in $\omega$
(see Figure \ref{f_many_chords}). Here, $\iota\in\{0,1,2,\ldots\}$
is a random variable which denotes the number of induced chords
in $\omega$ so that $\iota=0$ when the chord $(\xi'_0,\xi'_1)$
has no intersection with $\omega$.

%f5 ###
%

The generalization of Theorem 2.8 of \cite{CPSV1} that we
want to obtain is the following fact:
\begin{prop}
\label{p_induced_chords}
For any bounded function $f\dvtx{\mathfrak M}\mapsto{\mathbb R}$ we have
%e63 ###
%
\begin{eqnarray}\label{eq_induced_chords}\quad
&&\bigl\langle{\mathtt E}_{\omega}\bigl(f(\omega, {\mathcal U}\xi_0,
\xi_1)\mid\xi_0=(0,u)
\bigr) \bigr\rangle_{{\mathbb Q}}\nonumber\\
&&\qquad= \frac{{\mathcal Z}'}{{\mathcal Z}}\times\frac{1}{{\mathcal Z}'}
\int_{\Omega^2}d\widetilde{{\mathbb P}}(\omega,\omega')
\int_\Lambda d\mu_0^{\omega'} (u)\,\kappa_{0,u}^{-1} \\
&&\qquad\quad{} \times
{\mathtt E}_{\omega,\omega'} \Biggl(\sum_{k=1}^\iota
f \bigl(\theta_{\xi_0^{(k)}\cdot e}\omega, {\mathcal U}\xi_0^{(k)},
\xi_1^{(k)} - \bigl(\xi_0^{(k)}\cdot e\bigr)e \bigr)
\biggm| \xi'_0=(0,u) \Biggr).\nonumber
\end{eqnarray}
\end{prop}
\begin{pf}
We keep the notation from the proof of
Proposition \ref{p_cosine_int} (with the obvious modifications
in the case when $\omega'$ is considered instead of $\omega$).
Without restriction of generality, we suppose that $f$
is nonnegative.
First, analogously to (\ref{conta_inv*}),
we obtain that the right-hand side of (\ref{eq_induced_chords})
may be rewritten as
%e64 ###
%
\begin{eqnarray} \label{mnogo_chords2}
&&\frac{1}{2s{\mathcal Z}}
\int_{\Omega^2}d\widetilde{{\mathbb P}}(\omega,\omega')
\int_{{\tilde F}^{\omega'}_s}d\nu^{\omega'} (x)
\int_{\partial\omega'}d\nu^{\omega'} (y)\,K(x,y)\nonumber\\[-8pt]\\[-8pt]
&&\qquad{} \times\sum_{k=1}^{\iota(x,y)}
f \bigl(\theta_{x^{(k)}\cdot e}\omega, {\mathcal U}x^{(k)},
y^{(k)} - \bigl(x^{(k)}\cdot e\bigr)e \bigr) =: T_1,\nonumber
\end{eqnarray}
where $(x^{(1)},y^{(1)}),\ldots, (x^{(\iota(x,y))},
y^{(\iota(x,y))})$
are the chords induced in $\omega$ by the chord $(x,y)$
in $\omega'$.

Let us denote ${\tilde F}^\omega_{h,h'} =
\{x\in\partial\omega\dvtx x\cdot e\in[h,h']\}$
(so that ${\tilde F}^\omega_{s}={\tilde F}^\omega_{-s,s}$).
Define ${\hat\iota}_n(x,y)$ as the number of intersections
of the chord $(x,y)$ with ${\tilde F}^\omega_{s-n,s-n+1}$.
To proceed, we need the following fact: let $A$ be a
subset of $\partial{\mathcal D}_s^\omega$ and $x\in\partial
{\mathcal D}_s^{\omega'}$.
Then we have
\[
\mathtt{P}_{\omega'}[\mbox{random chord beginning at $x$
crosses $A$}]
\leq\nu^{\omega}(A)\sup_{y\in A}{\hat K}(x,y).
\]
Also, by decomposing $A$ into pieces that may have at most one
intersection with the chord starting from $x$ and using the above
inequality, we obtain
%e65 ###
%
\begin{eqnarray}\label{cross_a_piece}
&&\mathtt{E}_{\omega'}[\mbox{number of intersections of the
random chord from $x$
with $A$}]\nonumber\\[-8pt]\\[-8pt]
&&\qquad \leq\nu^{\omega}(A)\sup_{y\in A}{\hat K}(x,y).\nonumber
\end{eqnarray}
Using Condition \ref{ConditionL}
one obtains that $\nu^{\omega}({\tilde F}^\omega_{s-n,s-n+1})$ is bounded
from above by a constant [see the argument before (\ref{IP<IQ})].
From (\ref{def_K}) we know that $K(z,z')\leq
\frac{\gamma_d}{\|z-z'\|}$,
so for any $x\in\{s\}\times\omega'_s$ it is straightforward
to obtain that
%e66 ###
%
\begin{equation}
\label{number_inters}
\int_{\partial{\mathcal D}_s^{\omega'}}{\hat\iota}_n(x,y)K(x,y)\,
d\nu^{\omega'}(y) \leq
\frac{C_6\nu^{\omega}({\tilde F}^\omega_{s-n,s-n+1})}{n}
\leq\frac{C_7}{n}.
\end{equation}
Suppose, without restriction of generality, that $s$ is an integer
number.
Since $\iota(x,y) \leq1+\sum_{n=1}^{2s}{\hat\iota}_n(x,y)$,
we obtain that
%e67 ###
%
\begin{eqnarray}\label{laterals_not_count}
&&\frac{1}{s}\int_{\Omega^2}d\widetilde{{\mathbb P}}(\omega,\omega')
\int_{\{s\}\times\omega'_s}d\nu^{\omega'} (x)
\int_{\partial{\mathcal D}_s^{\omega'}} d\nu^{\omega'} (y)\,
{\hat K}(x,y)\iota(x,y)
\nonumber\\
&&\qquad\leq\frac{1}{s} \Biggl(1+\sum_{n=1}^{2s}\frac{C_7}{n} \Biggr)
\\
&&\qquad\leq\frac{C_8\ln s}{s}\nonumber
\end{eqnarray}
and the same bound also holds if we change $\{s\}\times\omega'_s$
to $\{-s\}\times\omega'_{-s}$ in the second integral above.

Note that, by the ergodic theorem, we have that
\[
\frac{\nu^{\omega}(\partial{\mathcal D}_s^{\omega})}{2s}\to
{\mathcal Z},\qquad
\frac{\nu^{\omega'} (\partial{\mathcal D}_s^{\omega'})}{2s}\to
{\mathcal Z}'\qquad
\mbox{as $s\to\infty$, $ \widetilde{{\mathbb P}}$-a.s.}
\]
Then, analogously to (\ref{conta_inv**}),
using (\ref{laterals_not_count})
together with the fact that $f$ is a bounded function,
we obtain that for any $\varepsilon>0$ there exists $s_0$ such that
for all $s\geq s_0$ [recall (\ref{mnogo_chords2})],
%e68 ###
%
\begin{equation}
\label{difference_Ts}
T_2-T_1<\varepsilon,
\end{equation}
where
%e69 ###
%
\begin{eqnarray} \label{mnogo_chords3}
&&T_2:=\frac{1}{\nu^{\omega'} (\partial{\mathcal
D}_s^{\omega'})}
\int_{(\partial{\mathcal D}_s^{\omega'})^2}d\nu^{\omega'} (x)\,
d\nu^{\omega'} (y)\,{\hat K}(x,y)\nonumber\\[-8pt]\\[-8pt]
&&\qquad{} \times\sum_{k=1}^{\iota(x,y)}
f \bigl(\theta_{x^{(k)}\cdot e}\omega, {\mathcal U}x^{(k)},
y^{(k)} - \bigl(x^{(k)}\cdot e\bigr)e \bigr).\nonumber
\end{eqnarray}
Then, by Theorem 2.8 of \cite{CPSV1}, we have
%e70 ###
%
\begin{equation}
\label{apply_Th.2.8}\hspace*{28pt}
T_2 = \frac{1}{\nu^{\omega}(\partial{\mathcal D}_s^{\omega})}
\int_{(\partial{\mathcal D}_s^{\omega})^2}d\nu^{\omega}(x)\, d\nu
^{\omega}(y)\,
{\hat K}(x,y) f \bigl(\theta_{x\cdot e}\omega, {\mathcal U}x,
y - (x\cdot e)e \bigr).
\end{equation}
Again, it is straightforward to obtain that for any $\varepsilon>0$
there exists $s_0$ such that for all $s\geq s_0$,
%e71 ###
%
\begin{eqnarray}\label{difference_others}\hspace*{22pt}
&& \biggl| \frac{1}{\nu^{\omega}(\partial{\mathcal D}_s^{\omega})}
\int_{(\partial{\mathcal D}_s^{\omega})^2}d\nu^{\omega}(x)\, d\nu
^{\omega}(y)\,
{\hat K}(x,y) f \bigl(\theta_{x\cdot e}\omega, {\mathcal U}x,
y - (x\cdot e)e \bigr)\nonumber\\
&&\qquad{}- \frac{1}{2s{\mathcal Z}}
\int_{{\tilde F}^\omega_s}d\nu^{\omega}(x)\\
&&\hspace*{85.7pt}{}\times\int_{\partial\omega} d\nu^{\omega}(y)\,
K(x,y) f \bigl(\theta_{x\cdot e}\omega, {\mathcal U}x,
y - (x\cdot e)e \bigr) \biggr| < \varepsilon.\nonumber
\end{eqnarray}
By the ergodic theorem, we have that ${\mathbb P}$-a.s.
\begin{eqnarray*}
&&\lim_{s\to\infty}\frac{1}{2s{\mathcal Z}}
\int_{{\tilde F}^\omega_s}d\nu^{\omega}(x)
\int_{\partial\omega} d\nu^{\omega}(y)
K(x,y) f \bigl(\theta_{x\cdot e}\omega, {\mathcal U}x,
y - (x\cdot e)e \bigr)\\
&&\qquad
= \bigl\langle{\mathtt E}_{\omega}\bigl(f(\omega, {\mathcal U}\xi_0, \xi
_1)\mid\xi_0=(0,u) \bigr)
\bigr\rangle_{{\mathbb Q}},
\end{eqnarray*}
so, using (\ref{difference_Ts}), (\ref{apply_Th.2.8})
and (\ref{difference_others}), we obtain,
abbreviating for a moment
\[
\mathfrak{A}:=\sum_{k=1}^\iota
f \bigl(\theta_{\xi_0^{(k)}\cdot e}\omega, {\mathcal U}\xi_0^{(k)},
\xi_1^{(k)} - \bigl(\xi_0^{(k)}\cdot e\bigr)e \bigr),
\]
that
%e72 ###
%
\begin{eqnarray}\label{ineq_induced_chords1}
&&\bigl\langle{\mathtt E}_{\omega}\bigl(f(\omega, {\mathcal U}\xi_0, \xi
_1)\mid\xi_0=(0,u) \bigr) \bigr\rangle
_{{\mathbb Q}}\nonumber\\[-8pt]\\[-8pt] %}\nonumber\\
&&\qquad\leq\frac{{\mathcal Z}'}{{\mathcal Z}}\times\frac{1}{{\mathcal Z}'}
\int_{\Omega^2}d\widetilde{{\mathbb P}}(\omega,\omega')
\int_\Lambda d\mu_0^{\omega'} (u)\,\kappa_{0,u}^{-1}
{\mathtt E}_{\omega,\omega'} \bigl(\mathfrak{A}
\mid\xi'_0=(0,u) \bigr).\nonumber
\end{eqnarray}
The other inequality is much easier to obtain.
Fix an arbitrary ${\tilde c}>0$, and consider
$\mathfrak{A} {\mathbb I}{\{\mathfrak{A}\leq{\tilde c}\}}$
instead of $\mathfrak{A}$. Since $\mathfrak{A} {\mathbb I}{\{
\mathfrak{A} \leq{\tilde c}\}}$ is bounded, we now have no
difficulties relating
the integrals on ${\tilde F}^{\omega'}_s\times\partial\omega'$
to the corresponding integrals on $(\partial{\mathcal D}_s^{\omega'})^2$.
In this way we obtain that for any ${\tilde c}$,
\begin{eqnarray*}
&&\bigl\langle{\mathtt E}_{\omega}\bigl(f(\omega, {\mathcal U}\xi_0, \xi
_1)\mid\xi_0=(0,u) \bigr) \bigr\rangle
_{{\mathbb Q}}\\ %}\nonumber\\
&&\qquad\geq\frac{{\mathcal Z}'}{{\mathcal Z}}\times\frac{1}{{\mathcal Z}'}
\int_{\Omega^2}d\widetilde{{\mathbb P}}(\omega,\omega')
\int_\Lambda d\mu_0^{\omega'} (u)\,\kappa_{0,u}^{-1}
{\mathtt E}_{\omega,\omega'}
\bigl(\mathfrak{A} {\mathbb I}{\{\mathfrak{A}\leq{\tilde c}\}}
\mid\xi'_0=(0,u) \bigr).
\end{eqnarray*}
We use now the monotone convergence theorem
and (\ref{ineq_induced_chords1}) to conclude the proof of
Proposition \ref{p_induced_chords}.
\end{pf}

Using Proposition \ref{p_induced_chords}, we are now
able to prove Proposition \ref{pr_suff_2nd_moment} for all
$d\geq3$.
\begin{pf*}{Proof of Proposition \protect\ref{pr_suff_2nd_moment}}
We apply Proposition \ref{p_induced_chords} with $\omega'$
being the straight cylinder, $\omega'={\mathbb R}\times\Lambda$.
For the random chord in a straight tube, using the fact that
the cosine density vanishes in the directions orthogonal to the
normal vector, we obtain that (for any starting point $\xi'_0$)
$\phi_0:= \mathtt{E}_{\omega'}((\xi'_1-\xi'_0)\cdot e)^2 < \infty$.

Now consider the function
$f_{\tilde c}(\omega,u,y)=(y\cdot e)^2{\mathbb I}{\{(y\cdot e)^2\leq
{\tilde c}\}}$.
Since
\[
\sum_{k=1}^\iota
f_{\tilde c} \bigl(\theta_{\xi_0^{(k)}\cdot e}\omega, {\mathcal U}\xi_0^{(k)},
\xi_1^{(k)} - \bigl(\xi_0^{(k)}\cdot e\bigr)e \bigr)
\leq\bigl((\xi'_1-\xi'_0)\cdot e\bigr)^2,
\]
we obtain that for any ${\tilde c}$,
\[
\bigl\langle{\mathtt E}_{\omega}\bigl(f_{\tilde c}(\omega, {\mathcal U}\xi_0,
\xi_1)\mid\xi
_0=(0,u) \bigr) \bigr\rangle_{{\mathbb Q}} \leq\phi_0.
\]
Using the monotone convergence theorem, we conclude
the proof of Proposition~\ref{pr_suff_2nd_moment}.
\end{pf*}

\textit{Remarks.}
(i) Observe from the definitions of $\phi_0$ above
and (\ref{def_Lambda}) of $\Lambda$ that $\phi_0({\widehat M})
\stackrel{\mathrm{def}}{=}
\phi_0={\widehat M}^2 \phi_0(1)$.
Then we have shown the universal bound
\[
\langle b \rangle_{{\mathbb Q}} \leq{\widehat M}^2 C(d)
\]
for all random tubes with a vertical section included in the
sphere $\Lambda$ of radius ${\widehat M}$ where
$C(d)=\phi_0(1)$ corresponds to the straight cylinder
with spherical section of radius ${\widehat M}=1$.

(ii) For $k\geq1$, denote by
\[
b^{(k)}(\omega,u) =
{\mathtt E}_{\omega}\bigl(|(\xi_1-\xi_0)\cdot e|^k \mid\xi_0=(0,u)\bigr)
\]
the $k$th absolute moment of the projection of the first jump to the
horizontal direction. Then, similarly to the proof of
Propositions \ref{pr_suff_2nd_moment}
and \ref{pr_infinite_d=2}, one can
obtain, for the $d$-dimensional model, that
$\langle b^{(d)} \rangle_{{\mathbb Q}}=+\infty$ in the case when the random
tube contains
a straight cylinder and that $\langle b^{(d-\delta)} \rangle
_{{\mathbb Q}
}<\infty$
for any $\delta>0$.

%s4.3 ###
\subsection{\texorpdfstring{Proof of Theorem
\protect\ref{t_q_invar_princ_cont}}{Proof of Theorem 2.2}}
\label{s_proof_inv_pr_cont}
We start by obtaining a formula for the mean length of the
first jump, at equilibrium.
\begin{lmm}
\label{l_chord_length}
We have
%e73 ###
%
\begin{equation}
\label{eq_chord_length}
\langle{\mathtt E}_{\omega}\|\xi_1-\xi_0\| \rangle_{{\mathbb Q}} =
\frac{\pi^{1/2}
\Gamma(({d+1})/{2})d}{\Gamma({d}/{2}+1)}\times
\frac{1}{{\mathcal Z}}\int_\Omega|\omega_0| \,d{\mathbb P}.
\end{equation}
\end{lmm}
\begin{pf}
Recall the notation ${\tilde F}^\omega_s$,
${\mathcal D}_s^\omega$, ${\hat K}(x,y)$
from the proof of Proposition~\ref{p_cosine_int}.
The strategy of the proof will be similar to that of the proof of
Proposition~\ref{p_induced_chords}.
Analogously to (\ref{conta_inv*}), we write
%e74 ###
%
\begin{eqnarray}\label{conta_dlina_chord}\qquad
\langle{\mathtt E}_{\omega}\|\xi_1-\xi_0\| \rangle_{{\mathbb Q}}
&=& \frac{1}{{\mathcal Z}}\int_\Omega d{\mathbb P}
\int_\Lambda
d\mu^{\omega}_0(u)\,\kappa_{0,u}^{-1}\nonumber\\
&&{}\times\int_{\partial\omega}d\nu
^{\omega}(y)\,
K ((0,u),y )\|(0,u)-y\|\nonumber\\
&=& \frac{1}{2s{\mathcal Z}}\int_\Omega d{\mathbb P}\int_{-s}^s ds
\int_\Lambda d\mu^{\omega}_s(u)\,\kappa_{s,u}^{-1}\\
&&{} \times\int_{\partial\omega}d\nu^{\omega}(y)\,
K ((s,u),y )
\|(s,u)-y\|\nonumber\\
&=& \frac{1}{2s{\mathcal Z}}\int_\Omega d{\mathbb P}\int_{{\tilde
F}^\omega_s}
d\nu^{\omega}(x)\int_{\partial\omega}d\nu^{\omega}(y) \|y-x\|K(x,y).\nonumber
\end{eqnarray}
By Theorem 2.6 of \cite{CPSV1}, we know that
%e75 ###
%
\begin{equation}
\label{finite_mean_chord}
\int_{(\partial{\mathcal D}_s^\omega)^2} d\nu^{\omega}(x)\, d\nu
^{\omega}(y)\,
{\hat K}(x,y) \|x-y\|
= \frac{\pi^{1/2}\Gamma(({d+1})/{2})d}
{\Gamma({d}/{2}+1)} |{\mathcal D}_s^\omega|.
\end{equation}
Denote by $D_\ell= \{-s\}\times\omega_{-s}$ and
$D_r = \{s\}\times\omega_{s}$ the left and
right vertical pieces of $\partial{\mathcal D}_s^\omega$, so that
$\partial{\mathcal D}_s^\omega= {\tilde F}^\omega_s\cup D_\ell\cup D_r$.
From (\ref{conta_dlina_chord})
we obtain [recall also that ${\hat K}(x,y)={\hat K}(y,x)$ for all
$x,y\in\partial{\mathcal D}_s^\omega$]
\begin{eqnarray*}
\langle{\mathtt E}_{\omega}\|\xi_1-\xi_0\| \rangle_{{\mathbb Q}}
&\geq&\frac{1}{2s{\mathcal Z}
}\int_\Omega d{\mathbb P}
\int_{({\tilde F}^\omega_s)^2}
d\nu^{\omega}(x)\, d\nu^{\omega}(y) \|y-x\|K(x,y) \\
&\geq&\frac{1}{2s{\mathcal Z}}\int_\Omega d{\mathbb P}
\biggl(\int_{(\partial{\mathcal D}_s^\omega)^2}
d\nu^{\omega}(x)\, d\nu^{\omega}(y) \|y-x\|{\hat K}(x,y) \\
&&\hspace*{54.8pt}{}
- 2 \int_{(D_\ell\cup D_r) \times
\partial{\mathcal D}_s^\omega}
d\nu^{\omega}(y) \|y-x\|{\hat K}(x,y) \biggr).
\end{eqnarray*}
Observe that [recall (\ref{def_Lambda})] for all
$x,y\in{\mathbb R}\times\Lambda$ it holds that
$\|x-y\|\leq|(x-y)\cdot e|+2{\widehat M}$.
So by (\ref{bound_d>=4}), there exists $C_1>0$ such that
for all $s$ we have
\[
\int_{(D_\ell\cup D_r) \times\partial{\mathcal D}_s^\omega}
d\nu^{\omega}(y) \|y-x\|{\hat K}(x,y) \leq C_1 \ln s + 2{\widehat M},
\]
and, using (\ref{finite_mean_chord}), we obtain
%e76 ###
%
\begin{eqnarray}
\label{oc_snizu_mean_chord_lim}
&&\langle{\mathtt E}_{\omega}\|\xi_1-\xi_0\| \rangle_{{\mathbb Q}}
\nonumber\\[-8pt]\\[-8pt]
&&\qquad\geq\frac{1}{2s{\mathcal Z}
}\int_\Omega d{\mathbb P}
\biggl(\frac{\pi^{1/2}\Gamma(({d+1})/{2})d}
{\Gamma({d}/{2}+1)}|{\mathcal D}_s^\omega|-C_1 \ln s
- 2{\widehat M} \biggr)\,d{\mathbb P}.\nonumber
\end{eqnarray}
Since, by the ergodic theorem,
\[
\frac{1}{2s}|{\mathcal D}_s^\omega| \to\int_\Omega|\omega_0|\,
d{\mathbb P}\qquad
\mbox{a.s., as $s\to\infty$},
\]
and sending $s$ to $\infty$
we obtain from (\ref{oc_snizu_mean_chord_lim}) that
%e77 ###
%
\begin{equation}
\label{oc_snizu_mean_chord}
\langle{\mathtt E}_{\omega}\|\xi_1-\xi_0\| \rangle_{{\mathbb Q}}
\geq\frac{\pi^{1/2}
\Gamma(({d+1})/{2})d}{\Gamma({d}/{2}+1)}\times
\frac{1}{{\mathcal Z}}\int_\Omega|\omega_0| \,d{\mathbb P}.
\end{equation}

Now, fix ${\tilde c}>0$ and write analogously
to (\ref{conta_dlina_chord})
\begin{eqnarray*}
&&\langle{\mathtt E}_{\omega}\|\xi_1-\xi_0\|{\mathbb I}{\{\|
\xi_1-\xi _0\| \leq {\tilde c}\}} \rangle_{{\mathbb Q}} \\
&&\qquad= \frac{1}{2s{\mathcal Z}}\int_\Omega d{\mathbb P}\int_{{\tilde
F}^\omega_s}
d\nu^{\omega}(x)\int_{\partial\omega}d\nu^{\omega}(y) \|y-x\|
{\mathbb I}{\{\|y-x\| \leq{\tilde c}\}} K(x,y).
\end{eqnarray*}
In this situation
$\|\xi_1-\xi_0\|{\mathbb I}{\{\|\xi_1-\xi_0\|\leq{\tilde c}\}}$ is bounded.
So, analogously to the proof of Proposition \ref{p_cosine_int} and
again using (\ref{finite_mean_chord}), by a coupling argument it is
elementary to obtain that for any ${\tilde c}$,
\begin{eqnarray*}
&&\langle{\mathtt E}_{\omega}\|\xi_1-\xi_0\|{\mathbb I}{\{\|\xi
_1-\xi_0\|\leq {\tilde c}\}}
\rangle_{{\mathbb Q}} \\
&&\qquad\leq
\frac{\pi^{1/2}\Gamma(({d+1})/{2})d}{\Gamma({d}/{2}+1)}
\times
\frac{1}{{\mathcal Z}}\int_\Omega|\omega_0| \,d{\mathbb P}.
\end{eqnarray*}
Using the monotone convergence theorem
and (\ref{oc_snizu_mean_chord}),
we conclude the proof of Lemma \ref{l_chord_length}.
\end{pf}

With Lemma \ref{l_chord_length} at hand,
we are now ready to prove Theorem \ref{t_q_invar_princ_cont}.
\begin{pf*}{Proof of Theorem \protect\ref{t_q_invar_princ_cont}}
Define $n(t)=\max\{n\dvtx\tau_n\leq t\}$.
We have
\[
t^{-1/2} X_t\cdot e = t^{-1/2} \xi_{n(t)}\cdot e
+ t^{-1/2}\bigl(X_t-\xi_{n(t)}\bigr)\cdot e.
\]
Let us prove first that the second term goes to $0$.
Indeed, by definition of the continuous-time process we have
%e78 ###
%
\begin{equation}
\label{eq_cont_clt}
t^{-1} \bigl(\bigl(X_t-\xi_{n(t)}\bigr)\cdot e \bigr)^2
\leq\frac{1}{n(t)} \bigl(\bigl(\xi_{n(t)+1}-\xi_{n(t)}\bigr)\cdot e
\bigr)^2.
\end{equation}
But then from the stationarity of $\xi$ we obtain that
\[
n^{-1} {\mathtt E}_{\omega}\bigl(\bigl(\xi_{n+1}-\xi_n\bigr)\cdot e \bigr)^2 \to0
\]
as $n\to\infty$ for ${\mathbb P}$-almost all $\omega$
[this is analogous to the derivation of (\ref{2nd_to_0})
in the proof of Theorem \ref{t_q_invar_princ}].

Now, the first term in the right-hand side of (\ref{eq_cont_clt})
equals
\[
\biggl(\frac{n(t)}{t} \biggr)^{1/2} \times
\frac{1}{n(t)^{1/2}}\xi_{n(t)}\cdot e.
\]
For the second term in the above product, apply
Theorem \ref{t_q_invar_princ}.
As for the first term, since
\[
\frac{n(t)}{\tau_{n(t)+1}} \leq\frac{n(t)}{t}
\leq\frac{n(t)}{\tau_{n(t)}},
\]
by the ergodic theorem and Lemma \ref{l_chord_length}
we have, almost surely,
\begin{eqnarray*}
\lim_{t\to\infty} \frac{n(t)}{t} &=&
\lim_{n\to\infty}\frac{n}{\tau_n}\\
& = &(\langle{\mathtt E}_{\omega}\|\xi_1-\xi_0\| \rangle_{{\mathbb
Q}} )^{-1}\\
& = & \frac{\Gamma({d}/{2}+1){\mathcal Z}}
{\pi^{1/2}\Gamma(({d+1})/{2})d}
\biggl(\int_\Omega|\omega_0|\, d{\mathbb P}\biggr)^{-1}.
\end{eqnarray*}
This concludes the proof of Theorem
\ref{t_q_invar_princ_cont}.\vspace*{-15pt}
\end{pf*}

\begin{appendix}\label{app}
\section*{Appendix}
\label{s_vertical_walls}
In this section we discuss the case when the map
$\alpha\mapsto\omega_\alpha$
is not necessarily continuous which corresponds to the case when
the random tube can have vertical walls.
The proofs contained here are given in a rather
sketchy way without going into much detail.

Define
\[
\varpi_\alpha= \{u\in\Lambda\dvtx (\alpha,u)\in\partial\omega\}
\]
to be the section of the boundary
by the hyperplane where the first coordinate is equal to $\alpha$.
Then let
\[
{\mathcal S}^{(0)} = \{\alpha\in{\mathbb R}\dvtx |\varpi_\alpha|\geq1\}
\]
and, for $n\geq1$,
\[
{\mathcal S}^{(n)} = \biggl\{\alpha\in{\mathbb R}\dvtx
|\varpi_\alpha|\in\biggl[\frac{1}{n+1},\frac{1}{n} \biggr) \biggr\}.
\]

Besides Condition \ref{ConditionR}, we have to assume something more.
Namely, we assume that for ${\mathbb P}$-almost all $\omega$,
$\nu^{\omega}$-almost all $(\alpha,u)\in\partial\omega$
are such that either $|\varpi_\alpha|>0$
(so that $\alpha\in{\mathcal S}^{(n)}$ for some $n$), or
$(\alpha,u)\in{\mathcal R}_{\omega}$ (recall the definition of
${\mathcal R}_{\omega}$ from
Section~\ref{s_notations}).

Also, we modify the definition of the measure $\mu^{\omega}_\alpha$
in the following way:
it is defined as in Section \ref{s_notations} when
$|\varpi_\alpha|=0$, and we put $\mu^{\omega}_\alpha\equiv0$
when $|\varpi_\alpha|>0$.

Observe that, for any $n\geq0$, ${\mathcal S}^{(n)}$ is a stationary point
process.
Note that, in contrast, the set $\bigcup_{n \geq0} {\mathcal S}^{(n)}$ may
\textit{not} be locally finite, which is the reason why we need to
introduce a sequence ${\mathcal S}^{(n)}$. Let ${\mathbb P}^{(n)}$
be the Palm version of ${\mathbb P}$ with respect to ${\mathcal
S}^{(n)}$, that is,
intuitively it is ${\mathbb P}$ ``conditioned on having a point of
${\mathcal S}^{(n)}$
at the origin.'' Observe that ${\mathbb P}^{(n)}$ is singular with respect
to ${\mathbb P}$, since, obviously,
\[
%% \label{wall_in_0}
{\mathbb P}[|\varpi_0|>0] = 0.
\]

Now, define (after checking that the two limits below exist
${\mathbb P}$-a.s.)
\begin{eqnarray*}
q_0 &=& \biggl(\int_\Omega|\varpi_0| \,d{\mathbb P}^{(0)}
\biggr)^{-1}\\
&&{}\times
\lim_{a\to\infty}\frac{\nu^{\omega}(\{x\in\partial\omega\dvtx
x\cdot e\in[0,a], |\varpi_{x\cdot e}|\geq1\} )}
{a},\\
q_n &=& \biggl(\int_\Omega|\varpi_0| \,d{\mathbb P}^{(n)} \biggr)^{-1}
\\
&&{}\times
\lim_{a\to\infty}\frac{\nu^{\omega}(\{x\in\partial\omega\dvtx
x\cdot e\in[0,a],
|\varpi_{x\cdot e}|\in[{1}/({n+1}),
{1}/{n})\} )}{a}
\end{eqnarray*}
for $n\geq1$.
Then, we define the measure ${\mathbb Q}$ which is the reversible measure
for the environment seen from the particle
%e79 ###
%
\begin{equation}\quad
\label{df_IQ_vert}
d{\mathbb Q}(\omega,u) = \frac{1}{{\mathcal Z}} \Biggl[\kappa_{0,u}^{-1}
\,d\mu^{\omega}_0(u)\,
d{\mathbb P}(\omega)
+ \sum_{n=0}^\infty q_n {\mathbb I}{\{u\in\varpi_0\}} \,du\,
d{\mathbb P}^{(n)}(\omega) \Biggr],
\end{equation}
where ${\mathcal Z}$ is the normalizing constant; as we will see below,
${\mathcal Z}$ still can be defined directly through the limit
\[
%% \label{df_norm_ZZ}
{\mathcal Z}= \lim_{a\to\infty}\frac{\nu^{\omega}(\{x\in\partial
\omega\dvtx
x\cdot e\in[0,a]\} )}{a}.
\]

The scalar product is now defined by
\begin{eqnarray*}
%% \label{df_scalar_vert}
\langle f, g \rangle_{{\mathbb Q}} &=& \frac{1}{{\mathcal Z}} \Biggl[\int
_\Omega d{\mathbb P}
\int_{\Lambda}
d\mu^{\omega}_0(u)\,
\kappa_{0,u}^{-1}f(\omega,u)g(\omega,u) \\
&&\hspace*{15.4pt} {}
+ \sum_{n=0}^\infty q_n \int_\Omega d{\mathbb P}^{(n)}
\int_{\varpi_0} du\, f(\omega,u)g(\omega,u) \Biggr].
\end{eqnarray*}
Now we need a slightly different definition for the transition
density: define $\bar{K}(x,y)$ by formula (\ref{def_K})
but without the indicator functions that $|{\mathbf n}_{\omega
}(x)\cdot e|\neq1$ and
$|{\mathbf n}_{\omega}(y)\cdot e|\neq1$. Also, the transition
operator $G$ can be
written in the following way:
\begin{eqnarray*}
Gf (\omega,u) &=& {\mathtt E}_{\omega}\bigl(f(\theta_{\xi_1\cdot e}\omega
,{\mathcal U}\xi_1)
\mid\xi_0 = (0,u)\bigr)\\
&=& \int_{\partial\omega} \bar{K} ((0,u),x )
f(\theta_{x\cdot e}\omega,{\mathcal U}x)\, d\nu^{\omega}(x)
\\
&=& \int_{-\infty}^{+\infty} d\alpha\int_{\Lambda}
d\mu^{\omega}_\alpha(v)\,
\kappa_{\alpha,v}^{-1}
f(\theta_\alpha\omega,v) \bar{K} ((0,u),(\alpha,v) )
\\
&&{} + \sum_{n=0}^\infty\sum_{\alpha\in{\mathcal S}^{(n)}}
\int_{\varpi_\alpha} dv\,
f(\theta_\alpha\omega,v) \bar{K} ((0,u),(\alpha,v) ).
%% \label{def_trans_operator_vert}
\end{eqnarray*}

Now, we have to prove the reversibility of $G$ with
respect to ${\mathbb Q}$.
The direct method adopted in the proof of
Lemma \ref{l_revers_discr} now seems to be to cumbersome to apply,
so we use another approach by taking advantage of stationarity.
For two bounded functions $f,g$, consider the quantity
\begin{eqnarray*}
{\mathfrak A}(a) &=& \frac{1}{{\mathcal Z}a}\int_{\{x\in\partial
\omega\dvtx
x\cdot e\in[0,a]\}^2}
d\nu^{\omega}(x)\,d\nu^{\omega}(y)\, \bar{K}(x,y)\\
&&\hspace*{91.3pt}{}\times f(\theta_{x\cdot
e}\omega,{\mathcal U}x)
g(\theta_{y\cdot e}\omega,{\mathcal U}y).
\end{eqnarray*}
Using (\ref{chord_not_too_long}),
it is elementary to obtain that (assuming for now that the
limit exists ${\mathbb P}$-a.s.)
\begin{eqnarray*}
\lim_{a\to\infty} {\mathfrak A}(a) &=& \lim_{a\to\infty}
\frac{1}{{\mathcal Z}a}\int_{\{x\in\partial\omega\dvtx x\cdot e\in
[0,a]\}}
d\nu^{\omega}(x)\,
f(\theta_{x\cdot e}\omega,{\mathcal U}x) \\
&&\hspace*{21.2pt}{} \times\int_{\partial\omega}d\nu^{\omega}(y)\,
\bar{K}(x,y)g(\theta_{y\cdot e}\omega,{\mathcal U}y)
\\
&=& \lim_{a\to\infty} \frac{1}{{\mathcal Z}a}\int_{\{x\in\partial
\omega\dvtx
x\cdot e\in[0,a]\}}d\nu^{\omega}(x)\,
f(\theta_{x\cdot e}\omega,{\mathcal U}x)\\
&&\hspace*{110.3pt}{}\times Gg(\theta_{x\cdot e}\omega
,{\mathcal U}x).
%% \label{navesim_hvosty}
\end{eqnarray*}
Then we write, using the ergodic theorem,
\begin{eqnarray*}
&&\lim_{a\to\infty} \frac{1}{a} \int_0^a d\alpha
\int_\Lambda d\mu^{\omega}_\alpha(u)\, \kappa^{-1}_{\alpha,u}
f(\theta_\alpha\omega,u)Gg(\theta_\alpha\omega,u)\\
&&\qquad=\int_\Omega d{\mathbb P}\int_\Lambda d\mu^{\omega}_0(u)\, \kappa^{-1}_{0,u}
f(\omega,u)Gg(\omega,u),\qquad {\mathbb P}\mbox{-a.s.}
%% \label{shod_nepr}
\end{eqnarray*}
Again, by the ergodic theorem, we have
\[
\lim_{a\to\infty} \frac{|{\mathcal S}^{(m)}\cap[0,a]|}{a} = q_m,\qquad
{\mathbb P}\mbox{-a.s.}
\]
so that we can write
\begin{eqnarray*}
&&\lim_{a\to\infty} \frac{1}{a}
\sum_{\alpha\in{\mathcal S}^{(m)}\cap[0,a]} \int_{\varpi_\alpha}du\,
f(\theta_\alpha\omega,u)Gg(\theta_\alpha\omega,u)\\
&&\qquad= q_m \int_\Omega d{\mathbb P}^{(m)}\int_{\varpi_0} du\,
f(\omega,u)Gg(\omega,u),\qquad {\mathbb P}\mbox{-a.s.}
%% \label{shod_diskr}
\end{eqnarray*}
Thus we have for ${\mathbb P}$-almost all environments
\begin{eqnarray*}
\lim_{a\to\infty} {\mathfrak A}(a) &=&\lim_{a\to\infty}
\frac{1}{{\mathcal Z}a}
\Biggl[\int_0^a d\alpha\int_\Lambda
d\mu^{\omega}_\alpha(u)\, \kappa^{-1}_{\alpha,u}
f(\theta_\alpha\omega,u)Gg(\theta_\alpha\omega,u) \\
&&\hspace*{45.3pt}{} + \sum_{m=0}^\infty
\sum_{\alpha\in{\mathcal S}^{(m)}\cap[0,a]} \int_{\varpi_\alpha}du\,
f(\theta_\alpha\omega,u)Gg(\theta_\alpha\omega,u) \Biggr] \\
& =& \frac{1}{{\mathcal Z}} \Biggl[\int_\Omega d{\mathbb P}\int_\Lambda
d\mu^{\omega}_0(u)\,
\kappa^{-1}_{0,u}
f(\omega,u)Gg(\omega,u) \\
&&\hspace*{15pt}{} + \sum_{m=0}^\infty q_m \int_\Omega d{\mathbb P}^{(m)}
\int_{\varpi_0} du\,
f(\omega,u)Gg(\omega,u) \Biggr]\\
&=& \langle f, Gg \rangle_{{\mathbb Q}}.
\end{eqnarray*}
By symmetry, in the same way one proves that $\lim_{a\to\infty}
{\mathfrak A}(a) = \langle g, Gf \rangle_{{\mathbb Q}}$, so $G$ is still
reversible with
respect to ${\mathbb Q}$.

Now the crucial observation is that formula (\ref{IP<IQ}) is
still valid even in the case when ${\mathbb Q}$ is defined
by (\ref{df_IQ_vert}), since we still have, for any $f\geq0$,
\[
\langle f \rangle_{{\mathbb Q}} \geq\frac{1}{{\mathcal Z}}\int
_{\Omega}d{\mathbb P}
\int_{\Lambda} d\mu^{\omega}_0(u)\,
\kappa_{0,u}^{-1}
f(\omega,u),
\]
so one can see that the whole argument goes through in this
general case as well.
However, we decided to write the proofs for the case of random
tube without vertical walls to avoid complicating the calculations
which are already quite involved. Here is the (incomplete) list of
places that would require modifications (and strongly complicate the
exposition):
\begin{itemize}
\item the display after (\ref{calc_g_nabla_f})
[part of the proof of (\ref{gradient=-diverg})];
\item the proof of (\ref{hatphi=b/2});
\item the proof of Proposition \ref{Prop_corrector};
\item the proof of (\ref{EQ_corr2<infty});
\item calculations (\ref{conta_inv*})
and (\ref{conta_dlina_chord}).
\end{itemize}
\end{appendix}

%%The work of F.C. was partially supported by CNRS (UMR 7599
%%``Probabilit{\'{e}}s et Mod{\`{e}}les Al{\'{e}}atoires'')
%%and ANR (grant POLINTBIO).
%%S.P. was partially supported by CNPq (300328/2005--2).
%Gunter M. Sch\"{u}tz thanks DFG (Schu 827/5--2, Priority programme SPP 1155) for
%financial support.
%%The work of M.V. was partially supported
%% by CNPq (304561/2006--1).
%Serguei Popov and Marina Vachkovskaia also thanks CNPq
%(471925/2006--3), FAPESP (04/07276--2) and CAPES/DAAD (Probral) for
%financial support.

% imsref loaded by lrinkeviciute, 2009-11-27 09:48:51
%

%
\printaddresses

\end{document}